\documentclass{article}

\usepackage{amssymb, latexsym}

\usepackage{graphicx}

\usepackage{amsmath}
\usepackage{color}

\newtheorem{theorem}{Theorem}[section]

\newtheorem{corollary}[theorem]{Corollary}

\newtheorem{lemma}[theorem]{Lemma}

\newtheorem{problem}[theorem]{Problem}

\newtheorem{proposition}[theorem]{Proposition}

\newtheorem{remark}[theorem]{Remark}

\begin{document}
\title{\large\bf  Harnack inequalities and $W$-entropy formula for Witten Laplacian on Riemannian manifolds with $K$-super Perelman Ricci flow}
\author{\ \ Songzi Li\thanks{Research partially supported by the China Scholarship Council.} , Xiang-Dong Li
\thanks{Research supported by NSFC No. 11371351, Key Laboratory RCSDS, CAS, No. 2008DP173182, and a
Hundred Talents Project of AMSS, CAS.} }

\maketitle

\begin{minipage}{120mm}
{\bf Abstract}.
In this paper, we prove logarithmic Sobolev inequalities and derive the Hamilton 
  Harnack inequality for the heat semigroup of the Witten Laplacian on complete Riemannian manifolds equipped with $K$-super Perelman
Ricci flow. We establish the $W$-entropy
formula for the heat equation of the Witten Laplacian and prove a rigidity theorem on complete Riemannian manifolds satisfying
the  $CD(K, m)$ condition,  and extend the $W$-entropy formula to time dependent Witten Laplacian on compact Riemannian manifolds with $(K, m)$-super Perelman Ricci flow, where  $K\in \mathbb{R}$ and $m\in [n, \infty]$ are two constants. Finally, we prove the 
Li-Yau and the Li-Yau-Hamilton  Harnack inequalities for positive solutions to the heat equation $\partial_t u=Lu$ associated to the time dependent Witten Laplacian on compact or complete manifolds equipped with variants of the $(K, m)$-super Ricci flow.

\end{minipage}

\tableofcontents

\section{Introduction}

\subsection{The differential Harnack inequality}

Differential Harnack inequality is an important tool in the study of heat equations and geometric flows on Riemannian manifolds. Let $M$ be an $n$ dimensional complete Riemannian manifold, $u$ be a
positive solution to the heat equation
\begin{eqnarray}
\partial_t u=\Delta u.\label{Heat1}
\end{eqnarray}
In their famous paper \cite{LY}, Li and Yau proved that, if $Ric\geq
-K$, where $K\geq 0$ is a positive constant, then for all
$\alpha>1$,
\begin{eqnarray}
{|\nabla u|^2\over u^2}-\alpha {\partial_t u\over u}\leq
{n\alpha^2\over 2t}+{n\alpha^2K\over \sqrt{2}(\alpha-1)}.\label{LYK}
\end{eqnarray}
In particular, if $Ric\geq 0$, then taking $\alpha\rightarrow 1$,
the Li-Yau differential Harnack inequality \cite{LY} holds
\begin{eqnarray}
{|\nabla u|^2\over u^2}-{\partial_t u\over u}\leq {n\over
2t}.\label{LY}
\end{eqnarray}

In \cite{H1}, Hamilton proved a dimension free Harnack inequality on
compact Riemannian manifolds with Ricci curvature bounded from
below. More precisely, if
\begin{eqnarray*}
Ric\geq -K,
\end{eqnarray*}
then, for any positive and bounded solution $u$ to the heat equation
$(\ref{Heat1})$, it holds
\begin{eqnarray}
{|\nabla u|^2\over u^2}\leq \left({1\over t}+2K\right)\log(A/u),\ \ \ \forall x\in M, t>0, \label{HamHar}
\end{eqnarray}
where
\begin{eqnarray*}
A:=\sup\limits\{u(t, x): x\in M, t\geq 0\}.
\end{eqnarray*}
Indeed, the same result holds on complete Riemannian manifolds with
Ricci curvature bounded from below. Under the same condition $Ric\geq -K$, Hamilton also
proved the following  Li-Yau type Harnack inequality  for any positive solution to the heat
equation $(\ref{Heat1})$
\begin{eqnarray}
{|\nabla u|^2\over u^2}-e^{2Kt}{\partial_t u\over u}\leq {n\over
2t}e^{4Kt}. \label{LYHHar}
\end{eqnarray}
In particular, when $K=0$, the above inequality reduces to the
Li-Yau Harnack inequality $(\ref{LY})$ on complete Riemannian
mnaifolds with non-negative Ricci curvature.

In this paper, we refer the inequality $(\ref{HamHar})$ Hamilton's (differential) Harnack inequality, and
refer the inequalities $(\ref{LY})$, $(\ref{LYK})$ and $(\ref{LYHHar})$ the Li-Yau-Hamiltom (differential) Harnack
inequality.

\subsection{The $W$-entropy formula}

Let $M$ be a closed manifold. In \cite{P1}, Perelman introduced the $\mathcal{F}$-entropy on the space of Riemannian metrics and smooth functions as follows
\begin{eqnarray*}
\mathcal{F}(g, f)=\int_M (R+|\nabla f|^2)e^{-f}dv,
\end{eqnarray*}
where $g\in \mathcal{M}=\{g:  {\rm Riemannian\ metric\ on}\ M\}$,
$f\in C^\infty(M)$, $R$ denotes the scalar curvature on $(M, g)$,
and $dv$ denotes the volume measure.  Under the constraint condition
that
\begin{eqnarray*}
dm=e^{-f}dv\end{eqnarray*}
is fixed, Perelman \cite{P1} proved that the gradient flow of $\mathcal{F}$ with respect to the standard $L^2$-metric on  $\mathcal{M}\times C^\infty(M)$ is given by the following modified Ricci flow for $g$ together with  the conjugate heat equation for $f$, i.e.,
\begin{eqnarray*}
\partial_t g&=&-2(Ric+\nabla^2 f),\\
\partial_t f&=&-\Delta f-R.
\end{eqnarray*}
Moreover, Perelman \cite{P1} introduced the
remarkable $W$-entropy as follows
\begin{eqnarray}
W(g, f, \tau)=\int_M \left[\tau(R+|\nabla
f|^2)+f-n\right]{e^{-f}\over
 (4\pi\tau)^{n/2}}dv,\label{entropy-1}
 \end{eqnarray}
where $\tau>0$, and $f\in C^\infty(M)$ satisfies the
following condition
$$
\int_M (4\pi\tau)^{-n/2}e^{-f}dv=1.$$  By \cite{P1}, it is known
that, if $(g(t), f(t), \tau(t))$ satisfies the evolution
equations
\begin{eqnarray}
\partial_t g=-2Ric, \ \ \partial_t f=-\Delta f+|\nabla
f|^2-R+{n\over 2\tau}, \ \ \partial_t \tau=-1, \label{r-c}
\end{eqnarray}
where the first one is Hamilton's Ricci flow, and the second one is
the corresponding conjugate heat equation, then the following
Perelman entropy formula holds
\begin{eqnarray}
{d\over dt}W(g, f, \tau)=2 \int_M \tau\left|Ric+\nabla^2
f-{g\over 2\tau}\right|^2{e^{-f}\over (4\pi \tau)^{n/2}}dv.\label{Entropy-P}
\end{eqnarray}
This implies that the $W$-entropy is increasing
in $\tau$ and the monotonicity is strict except that $M$ is a
shrinking Ricci soliton
\begin{eqnarray*}
Ric+\nabla^2 f={g\over 2\tau}.
\end{eqnarray*}
As an application of the above entropy formula, Perelman \cite{P1} derived the non local collapsing theorem for the Ricci flow,
 which plays an important r\^ole for ruling out cigars, the one part of the singularity
classification for the final resolution of the Poincar\'e conjecture and geometrization conjecture.

Since Perelman's preprint \cite{P1} was published on Arxiv in 2002, many
people have studied the $W$-like entropy for other geometric
flows on Riemannian manifolds \cite{N1, N2, Ec, LNVV, KN}. In
\cite{N1, N2}, Ni studied the $W$-entropy for the linear
heat equation on complete Riemannian manifolds.  More precisely, let
$(M, g)$ be an $n$-dimensional complete Riemannian manifold, let
$$u={e^{-f}\over (4\pi t)^{n/2}}$$
be a positive solution to the linear heat equation
\begin{eqnarray}
\partial_t u=\Delta u \label{heat-1}
\end{eqnarray}
with $\int_M u(x, 0)dv(x)=1$. The $W$-entropy  for the
linear heat equation $(\ref{heat-1})$ is defined by
\begin{eqnarray}
W(f, t)=\int_M \left[t |\nabla
f|^2+f-n\right]{e^{-f}\over
 (4\pi t)^{n/2}}dv.\label{entropy-2}
\end{eqnarray}
In \cite{N1}, Ni proved the following entropy formula
\begin{eqnarray}
{dW(f, t)\over dt}&=&-2\int_M t\left(
\left|\nabla^2 f-{g\over 2t}\right|^2+Ric(\nabla f, \nabla
f)\right){e^{-f}\over (4\pi t)^{n/2}}dv. \label{entropy-Ni}
\end{eqnarray}
This yields that  the $W$-entropy for the
linear heat equation $(\ref{heat-1})$ is decreasing on complete Riemannian
manifolds with non-negative Ricci curvature.

In \cite{Li12, Li13},  the second author of this paper introduced the $W$-entropy for the heat equation associated with the Witten Laplacian and proved the monotonicity and rigidity results on complete Riemannian manifolds with non-negative $m$-dimensional
Bakry-Emery Ricci curvature condition. More
precisely, let $(M, g)$ be a complete Riemannian manifold, $\phi\in C^2(M)$. Let
$$L=\Delta-\nabla\phi\cdot\nabla$$
be the Witten Laplacian on $(M, g)$ with respect to $\mu$, where
$$d\mu=e^{-\phi}dv.$$
For any $m\in [n, \infty)$, let
\begin{eqnarray*}
Ric_{m, n}(L)=Ric+\nabla^2\phi-{\nabla\phi\otimes\nabla\phi\over m-n},
\end{eqnarray*}
be the $m$-dimensional Bakry-Emery Ricci curvature of $L$, where $m=n$ if and only if $\phi$ is identically constant. Let
 $u={e^{-f}\over (4\pi t)^{m/2}}$ be the fundamental solution to
the heat equation
\begin{eqnarray*}
\partial_t u=Lu.
\end{eqnarray*}
Let
\begin{eqnarray}
H_m(u, t)=-\int_M u\log u d\mu-{m\over 2}(1+\log(4\pi t)).\label{Hm}
\end{eqnarray}
Define the $W$-entropy by the Boltzmann formula
\begin{eqnarray}
W_m(u, t)={d\over dt}(tH_m(u, t)).\label{Wm}
\end{eqnarray}
Then
\begin{eqnarray}
{d\over dt}H_m(u, t)=-\int_M \left(L\log u+{m\over
2t}\right)ud\mu,\label{LYHHm}
\end{eqnarray}
and
\begin{eqnarray}
W_m(u, t)=\int_M\left[t|\nabla \log u|^2+f-m\right]{e^{-f}\over
(4\pi t)^{m/2}}d\mu.\label{W-0}
\end{eqnarray}
Moreover, under the condition that  $(M, g)$ is a complete Riemannian manifold with bounded geometry condition, and $\phi\in C^4(M)$ with $\nabla\phi\in C_b^3(M)$, we have\footnote{In 2006,  the second author of this paper proved the $W$-entropy formula  $(\ref{W-1})$ for
 all positive solutions to the heat equation $\partial_t u=Lu$ on compact Riemannian manifolds with fixed metrics and potentials. See \cite{Li07}.}
\begin{eqnarray}
{d W_m(u, t)\over dt}&=&-2\int_M
t \left(\left|\nabla^2 f-{g\over 2t}\right|^2+Ric_{m,
n}(L)(\nabla f,
\nabla f)\right)ud\mu\nonumber\\
& & \hskip2cm  -{2\over m-n}\int_M t \left({\nabla\phi\cdot\nabla
f}+{m-n\over 2t}\right)^2ud\mu.\label{W-1}
\end{eqnarray}
In particular, if $(M, g, \phi)$ satisfies the bounded geometry condition and $Ric_{m, n}(L)\geq 0$, then the $W$-entropy is decreasing in time $t$, i.e.,
\begin{eqnarray*}
{dW_m(u, t)\over dt}\leq 0, \ \ \ \ \ \ \ \forall t\geq 0.
\end{eqnarray*}
Moreover,  ${dW_{m}(u, t)\over dt}=0$ holds at some $t=t_0>0$ if and only if $M$ is isomeric to $\mathbb{R}^n$, $m=n$, $\phi=C$ for some constant $C\in \mathbb{R}$, and $u(x, t)={e^{-{|x|^2\over 4t}}\over (4\pi t)^{n/2}}$ for all $x\in \mathbb{R}^n$ and $t>0$.

Here we say that $(M, g)$ satisfies the bounded geometry condition if the Riemannian curvature
 tensor ${\rm Riem}$ and its covariant derivatives $\nabla^k {\rm Riem}$ are uniformly bounded on
 $M$, $k=1, 2, 3$. The bounded geometry condition and the assumption $\phi\in C^4(M)$ with $\nabla\phi\in C^3_b(M)$ are only required in order to allow us to exchange the time derivatives and the integration of $u\log u$ on complete non-compact Riemannian manifolds.

In \cite{LL13}, when $m\in \mathbb{N}$,  we gave a direct proof of the $W$-entropy formula
$(\ref{W-1})$ for the Witten Laplacian by applying Ni's $W$-entropy formula
$(\ref{entropy-2})$  for the usual Laplacian to $M\times S^{m-n}$ equipped with a suitable warped product Riemannian metric, and
 gave a natural geometric interpretation for the third term in the $W$-entropy formula $(\ref{W-1})$
 for the Witten Laplacian. We have further proved
 the $W$-entropy formula for time dependent Witten Laplacian on Riemannian manifolds with time dependent metrics and potentials. In particular, if $d\mu=e^{-\phi}dv$ is fixed and if
\begin{eqnarray*}
{1\over 2}{\partial g\over \partial t}+Ric_{m, n}(L)\geq 0,
\end{eqnarray*}
then the $W$-entropy for the time dependent Witten Laplacian $L=\Delta_{g(t)}-\nabla_{g(t)}\phi(t)\cdot\nabla_{g(t)}$ is decreasing in time. For details, see \cite{LL13}.

\subsection{The purpose of our paper}

In \cite{Li05, BL, Li12}, the Li-Yau Harnack inequality $(\ref{LY})$
has been extended to positive solutions of the heat equation
associated to the Witten Laplacian on complete Riemannian manifolds
with non-negative
 $m$-dimensional Bakry-Emery Ricci curvature. More precisely, if
 $Ric_{m, n}(L)\geq 0$, or equivalently, if the  $CD(0, m)$ condition holds on $(M, g, \phi)$, then
\begin{eqnarray}
L\log u+{m\over 2t}\geq 0. \ \ \label{LYHm}
\end{eqnarray}
In \cite{Li13}, the Hamilton Harnack inequality has been also
extended to positive and bounded solutions to the heat equation of the Witten Laplacian on complete Riemannian manifolds
with the infinite dimensional Bakry-Emery Ricci curvature bounded
from below, or equivalently, on $(M, g, \phi)$ with the  $CD(-K, \infty)$ condition, i.e., $Ric(L)=Ric+\nabla^2\phi\geq -K$, where $K\in \mathbb{R}$ is a constant. As already mentioned  here,  the main ingredient in
these work is the introduction of the notion of finite or infinite
dimensional Bakry-Emery Ricci curvature, which plays an important
role to substitute that of Ricci curvature in the study of
geometric analysis for the usual Laplace-Beltrami operator on
Riemannian manifolds.

On the other hand, from $(\ref{Hm})$, $(\ref{Wm})$ and $(\ref{LYHHm})$, we can
see that there exists an essential and deep relationship among the $W$-entropy for the Witten Laplacian, the Gaussian heat kernel on $\mathbb{R}^m$ and  the Li-Yau Harnack inequality $(\ref{LYHm})$ on complete Riemannian manifolds satisfying the $CD(0,
m)$ condition.  Indeed,  the
$W$-entropy  for the heat equation $\partial_t u=Lu$ is defined by the Boltzmann formula
$(\ref{Wm})$ together with $(\ref{Hm})$, in which the quantity $H_{m}(u, t)=-\int_M u\log ud\mu-{m\over 2}(\log(4\pi t)+1)$ is defined as the difference between the Boltzmann-Shannon entropy $H(u)=-\int_M u\log ud\mu$ for the heat kernel measure $u(x, t)d\mu(x)$ of the Witten Laplacian on $M$ and the Boltzmann-Shannnon entropy $H(\overline{u})=-\int_{\mathbb{R}^m} \overline{u}\log \overline{u}dx={m\over 2}(\log(4\pi t)+1)$ of the Gaussian heat kernel measure $\overline{u}(x, t)dx$ on $\mathbb{R}^m$ with $\overline{u}(x, t)= {e^{-{\|x\|^2\over 4t}} \over (4\pi t)^{m/2}}$. By $(\ref{LYHHm})$, the time derivative of $H_{m}(u, t)$ is given by the integral of the Li-Yau Harnack quantity
$L\log u+{m\over 2t}$ with respect to the heat kernel measure $u(x, t)d\mu(x)$. The  $W$-entropy formula $(\ref{W-1})$ does not only imply the
monotonicity of the $W$-entropy on complete
Riemannian manifolds with the $CD(0, m)$ condition, but also allows us to prove a rigidity theorem which characterizes the unique equilibrium state of the $W$-entropy on the canonical ensemble of all complete Riemannian manifolds satisfying the $CD(0, m)$ condition. In \cite{LL13}, we defined the $W$-entropy  for the heat equation of the time dependent
Witten Laplacian by the same formulas  $(\ref{Hm})$ and $(\ref{Wm})$, and proved that
the $W$-entropy introduced in this way is monotonically decreasing in time on compact Riemannian
manifolds with ${1\over 2}{\partial g\over
\partial t}+Ric_{m, n}(L)\geq 0$.

Now it is very natural to raise the following problems: (1)  How to define the $W$-entropy functional for the heat equation associated with the Witten
Laplacian on complete Riemannian manifolds satisfying the $CD(K, \infty)$ condition or the $CD(K, m)$ condition for $K\in \mathbb{R}$ and $m\in [n, \infty)$? (2) Can we establish the monotonicity and rigidity  theorems for the $W$-entropy associated with the Witten
Laplacian on complete Riemannian manifolds satisfying general curvature-dimension condition? (3) What happens  on Riemannian manifolds with time dependent metrics and potentials?  (4) How to extend the Li-Yau and the Li-Yau-Hamilton 
Harnack inequalities to positive solutions to the heat equation associated to the time dependent Witten Laplacian on compact or complete manifolds with  time dependent metrics and potentials? Indeed, we have been asked these questions for many times by many people during the past years.

The purpose of this paper is to study these problems. We first prove the logarithmic Sobolev inequality and the reversal logarithmic Sobolev inequality and derive the Hamilton Harnack inequality  for the heat equation of the Witten Laplacian  on complete Riemannian manifolds with time
dependent metrics and potentials evolving along  $K$-super Perelman  Ricci flow. Then we introduce the $W$-entropy and prove the $W$-entropy
formula and rigidity theorem for the Witten Laplacian on complete
Riemannian manifolds with fixed metrics and potentials satisfying the
$CD(K, m)$ condition, for $K\in \mathbb{R}$ and $m\in
[n, \infty]$. Moreover, we extend the $W$-entropy results to compact Riemannian manifolds with time
dependent metrics and potentials evolving along $K$-super Perelman
Ricci flow. Finally, we prove the Li-Yau and Li-Yau-Hamilton Harnack inequalities  for positive solutions to the heat equation associated to the time dependent Witten Laplacian on compact or complete manifolds with $(K, m)$-super Ricci flow. 

\subsection{Statement of main results}

To state our results, let us first introduce some notations. Let $M$
be a complete Riemannian manifold with a fixed Riemnnian metric $g$,
$\phi\in C^2(M)$ and $d\mu=e^{-\phi}dv$, where $v$ is the Riemannian
volume measure on $(M, g)$. The Witten Laplacian  on $(M, g)$ with respect to the weighted volume measure $\mu$ or the potential function $\phi$ is defined  by
\begin{eqnarray*}
L =\Delta -\nabla \phi\cdot\nabla. \label{WL}
\end{eqnarray*}
For all $u, v\in C^\infty_0(M)$, the following integration by parts formula holds
\begin{eqnarray*}
\int_M \langle \nabla u, \nabla v\rangle d\mu=-\int_M Lu vd\mu=-\int_M uLvd\mu.
\end{eqnarray*} In \cite{BE}, Bakry and Emery proved that for all $u\in C_0^\infty(M)$,
\begin{eqnarray}
L|\nabla u|^2-2\langle \nabla u, \nabla L u\rangle=2|\nabla^2
u|^2+2Ric(L)(\nabla u, \nabla u), \label{BWF}
\end{eqnarray}
where
$$Ric(L)=Ric+\nabla^2\phi.$$
The formula $(\ref{BWF})$ can be viewed as a natural extension of the Bochner-Weitzenb\"ock formula. The quantity $Ric(L)=Ric+\nabla^2 \phi$, called the infinite dimensional Bakry-Emery Ricci curvature on the weighted Riemannian manifolds $(M, g, \phi)$. It plays as a good substitute of the Ricci curvature in many problems in comparison geometry and analysis on complete Riemannian manifolds with smooth weighted volume measures. See \cite{BE, BL, FLZ, FLL, Li05, Li12, Lot,  WW} and reference therein.

Following \cite{BE, Lot, Li05}, we introduce the $m$-dimensional Bakry-Emery Ricci curvature on $(M, g, \phi)$ by
\begin{eqnarray*}
Ric_{m, n}(L):=Ric+\nabla^2\phi-{\nabla\phi\otimes \nabla \phi\over m-n},
\end{eqnarray*}
where $m\geq n$ is a constant, and $m=n$ if and only if $\phi$ is a constant. When $m=\infty$, we have $Ric_{\infty, n}(L)=Ric(L)$.
Following \cite{BL}, we say that the Witten Laplacian $L$ satisfies the $CD(K, \infty)$ condition if $Ric(L)\geq K$, and $L$ satisfies the $CD(K, m)$ condition if $Ric_{m, n}(L)\geq K$. Recall that,
when $m\in \mathbb{N}$, the $m$-dimensional Bakry-Emery Ricci curvature $Ric_{m, n}(L)$ has a very natural geometric interpretation.
Indeed, consider the warped product metric on $M^n\times S^{m-n}$ defined by
\begin{eqnarray*}
\widetilde{g}=g_M\bigoplus e^{-{2\phi\over m-n}}g_{S^{m-n}}.\label{WPM}
\end{eqnarray*}
where $S^{m-n}$ is the unit sphere in $\mathbb{R}^{m-n+1}$ with the standard metric $g_{S^{m-n}}$. By a classical result in Riemaniann geometry, the quantity $Ric_{m, n}(L)$ is equal to the Ricci curvature of the above warped product metric $\widetilde g$ on $M^n\times S^{m-n}$ along the horizontal vector fields. See \cite{Lot, Li05, WW}. 

Let $(M, g(t), \phi(t), t\in [0, T])$ be a complete Riemannian manifold equipped with a family of time dependent Riemannian metrics $g(t)$ and potential functions $\phi(t)$, $t\in [0, T]$. In this paper, we call $(M, g(t), \phi(t), t\in [0, T])$ a $(K, m)$-super Ricci flow if the metric $g(t)$ and the potential function $\phi(t)$ satisfy the following inequality
\begin{eqnarray*}
{1\over 2}{\partial g\over \partial t}+Ric_{m, n}(L)\geq Kg.
\end{eqnarray*}
When $m=\infty$, i.e., if the metric $g(t)$ and the potential function $\phi(t)$ satisfy the following inequality
\begin{eqnarray*}
{1\over 2}{\partial g\over \partial t}+Ric(L)\geq Kg,
\end{eqnarray*}
we call $(M, g(t), \phi(t), t\in [0, T])$   a $K$-super Perelman Ricci flow.  As mentioned in Section $1.2$ and changing $f$ by $\phi$,  the modified Ricci flow ${1\over 2}{\partial g\over \partial t}+Ric(L)=0$ was introduced by 
Perelman \cite{P1} and can be regarded as the gradient flow of $\mathcal{F}(g, \phi)=\int_M (R+|\nabla \phi|^2)e^{-\phi}dv$ on $\mathcal{M}\times C^\infty(M)$ under the constraint condition that $dm=e^{-\phi}dv$ does not change in time. 
 
We now state the main results and describe the organisation of this paper. 

In Section $2$, we prove the following result  which describes the equivalence between the $K$-super Perelman Ricci flow property and the  logarithmic Sobolev inequality, the reversal  logarithmic Sobolev inequality for the heat semigroup associated
with the time dependent Witten Laplacian on manifolds with time dependent metrics and potentials.

\begin{theorem}\label{Thm-RLSI} Let $M$ be a complete Riemannian manifold
equipped with a family of time dependent metrics and
$C^2$-potentials $(g(t), \phi(t), t\in [0, T])$. Let
$L=\Delta_{g(t)}-\nabla_{g(t)}\phi(t)\cdot\nabla_{g(t)}$ be the time
dependent weighted Laplacian on $(M, g(t), \phi(t))$, $P_{s, t}f=u(\cdot, t)$ be a positive solution to the heat equation $\partial_t
u=L u$ with the initial condition $u(\cdot, s)=f$, where $0\leq s<t\leq T$, and $f$ is a
$C^1$ smooth and positive function on $M$. Then $(g(t), \phi(t), t\in [0,
T])$ satisfies a $K$-super Perelman Ricci flow equation
\begin{eqnarray}
\label{RFK} {1\over 2}{\partial g\over \partial t} +
Ric(L)\geq -K, \end{eqnarray}
 where $K\geq 0$ is a constant, if and only if for $0 \leq s < t \leq T$, the following logarithmic Sobolev inequality holds
\begin{eqnarray*}
P_{s, t}(f\log f)-P_{s, t}f\log P_{s,t}f\leq {e^{2K(t-s)}-1\over 2K}P_{s, t}\left({|\nabla f|^2\over
f}\right), \end{eqnarray*} or the reversal
logarithmic Sobolev inequality holds
\begin{eqnarray}
 {|\nabla P_{s,t} f|^2\over P_{s, t}
f}\leq {2K\over 1-e^{-2K(t-s)}}\left(P_{s, t}(f\log f)-P_{s,t}f\log P_{s, t}f\right).
\end{eqnarray}
\end{theorem}

As a corollary of the above theorem, we derive the following Hamilton Harnack inequality for the positive solution of the heat equation $\partial_t u=Lu$ on complete Riemannian manifolds with $K$-super Perelman Ricci flow. 

\begin{theorem}\label{Thm0} Let $M$ be a complete Riemannian manifold equipped with a family of time dependent metrics and $C^2$-potentials
$(g(t), \phi(t), t\in [0, T])$ satisfying a $K$-super Perelman
Ricci flow equation
\begin{eqnarray*}
{1\over 2}{\partial g\over \partial t}+Ric(L)\geq -K.
\end{eqnarray*}
where $K\geq 0$ is a constant independent of  $t\in [0, T]$. Let $u$ be a positive and bounded solution to the heat equation
\begin{eqnarray*}
\partial_t u=Lu,
\end{eqnarray*}
where
$$L=\Delta_{g(t)}-\nabla_{g(t)}\phi(t)\cdot\nabla_{g(t)}
$$
is the time dependent Witten Laplacian on $(M, g(t), \phi(t))$.
 Then for all $x\in M$ and $t>0$,
\begin{eqnarray}
{|\nabla u|^2\over u^2}\leq {2K\over 1-e^{-2Kt}}\log (A/u),\ \ \label{HH}
\end{eqnarray}
where
\begin{eqnarray*}
A:=\sup\limits\{u(t, x): x\in M, t\geq 0\}.
\end{eqnarray*}
In particular, the Hamilton Harnack inequality holds
\begin{eqnarray}
{|\nabla u|^2\over u^2}\leq \left({1\over
t}+2K\right)\log(A/u).\label{Ham}
\end{eqnarray}
In the case $K=0$, i.e.,
$(M, g(t), \phi(t), t\in [0, T])$ is a complete Riemannian manifold equipped with the super Perelman
Ricci flow
\begin{eqnarray*}
{1\over 2}{\partial g\over \partial t}+Ric(L)\geq 0,
\end{eqnarray*}
we have
\begin{eqnarray*}
{|\nabla u|^2\over u^2}\leq {1\over t}\log {A\over u}.\ \ \label{HH1}
\end{eqnarray*}

\end{theorem}

In particular, taking $\phi=0$, $m=n$ and $L=\Delta$, we have the Hamilton  Harnack inequality on complete Riemannian manifolds with $K$-super Ricci flow.

\begin{theorem}\label{Thm02}Let $M$ be a complete Riemannian manifold equipped with a family of Riemannian metrics $(g(t), t\in [0, T])$ evolving along a $K$-super Ricci flow
\begin{eqnarray*}
{1\over 2}{\partial g\over \partial t}+Ric\geq -K,
\end{eqnarray*}
where $K\geq 0$ be a constant independent of $t\in [0, T]$. Let $u$ be a positive and bounded solution to the heat equation
\begin{eqnarray*}
\partial_t u=\Delta u.
\end{eqnarray*}
Let $A:=\sup\limits\{u(t, x): x\in M, t\geq 0\}$. Then for all $x\in M$ and $t>0$,
\begin{eqnarray*}
{|\nabla u|^2\over u^2}\leq {2K\over 1-e^{-2Kt}}\log (A/u).\ \ \label{HH1}
\end{eqnarray*}
In particular, the Hamilton Harnack inequality holds
\begin{eqnarray*}
{|\nabla u|^2\over u^2}\leq \left({1\over
t}+2K\right)\log(A/u).\label{Ham1}
\end{eqnarray*}
In the case $K=0$, i.e.,
$(M, g(t), t\in [0, T])$ is a complete Riemannian manifold equipped with a super Ricci flow
\begin{eqnarray*}
{1\over 2}{\partial g\over \partial t}+Ric\geq 0,
\end{eqnarray*}
we have \footnote{In \cite{ZQ2}, Qi S. Zhang proved $(\ref{HH2})$ for the heat equation $\partial_t u=\Delta u$ on compact or complete Riemmanian manifolds equipped with the Ricci flow $\partial_t g=-2Ric$. }
\begin{eqnarray}
{|\nabla u|^2\over u^2}\leq {1\over t}\log {A\over u}.\ \ \label{HH2}
\end{eqnarray}

\end{theorem}

Moreover, we also prove the following theorem which extends the Li-Yau-Hamilton Harnack inequality $(\ref{LYHHar})$ to positive solutions of the heat equation $\partial_t u=Lu$ on complete Riemannian manifolds with
 fixed metrics and potentials satisfying the $CD(-K, m)$ condition.

\begin{theorem}\label{HLYH} Let $(M, g)$ be a complete Riemannian manifold with a $C^2$-potential $\phi$. Suppose that there exist some
constants $m\geq n$ and $K\geq 0$ such that
$$Ric_{m, n}(L)\geq -K.$$ Let $u$ be a
positive solution of the heat equation
\begin{eqnarray*}
\partial_t u=Lu.
\end{eqnarray*}
Then the Li-Yau-Hamilton Harnack inequality holds
\begin{eqnarray*}
{\partial_t u\over u}-e^{-2Kt}{|\nabla u|^2\over u^2}+e^{2Kt}{m\over
2t}\geq 0.
\end{eqnarray*}
In particular, if $K=0$, i.e., $Ric_{m, n}(L)\geq 0$, then the Li-Yau  Harnack
inequality holds
\begin{eqnarray*}
{\partial_t u\over u}-{|\nabla u|^2\over u^2}+{m\over 2t}\geq 0.
\end{eqnarray*}
\end{theorem}

As the corollaries of Theorem \ref{Thm0} and Theorem \ref{HLYH}, we have the following Harnack inequalities for positive solutions of the heat equation of the Witten Laplacian.

\begin{corollary}\label{cor1} Under the same condition as in Theorem \ref{Thm0}, for any $\delta>0$, and for all $x, y\in M$, $0<t<T$, we have
\begin{eqnarray*}
u(x, t)\leq u(y, t)^{1\over 1+\delta}A^{\delta\over 1+\delta}\exp\left\{{1+\delta^{-1}\over 4(1+\delta)}{2K\over 1-e^{-2Kt}}d^2(x, y)\right\}.
\end{eqnarray*}
\end{corollary}

\begin{corollary}\label{cor2} Let $M$ be a complete Riemannian manifold with $Ric_{m, n}(L)\geq -K$, $u$ be a positive solution to the heat equation $\partial_t u=Lu$. Then, for all $x, y\in M$, $0<\tau<T$, we have
\begin{eqnarray*}
u(x, \tau)\leq \left({T\over \tau}\right)^{m/2}u(y, T)\exp\left\{{1\over 4}e^{2K\tau}[1+2K(T-\tau){d^2(x, y)\over T-\tau}+{m\over 2}[e^{2KT}-e^{2K\tau}]\right\}.
\end{eqnarray*}

\end{corollary}

In Section $3$, we introduce the $W$-entropy and prove the $W$-entropy
formulas  for the heat equation of the Witten
Laplacian on complete Riemannian manifolds satisfying the $CD(K, m)$ condition, for $K\in \mathbb{R}$ and $m\in [n, \infty]$.  We also prove a rigidity theorem on complete Riemannian manifolds satisfying the $CD(K, m)$ condition for $K\in \mathbb{R}$ and $m\in [n, \infty)$. These extend the
$W$-entropy formula proved in \cite{Li12, Li13, LL13} for the Witten Laplacian
on complete Riemannian manifolds satisfying the $CD(0, m)$ condition, $m\in [n, \infty)$. We will also extend the $W$-entropy formula to time dependent Witten Laplacian on compact
Riemannian manifolds with $K$-super Perelman  Ricci flow.

\begin{theorem}\label{WCDK} Let $M$ be a complete Riemannian manifold with
bounded geometry condition, $\phi\in C^4(M)$ with $\nabla\phi\in
C_b^3(M)$. Suppose that $Ric+\nabla^2\phi\geq K$, where $K\in \mathbb{R}$ is
a constant. Let $u(\cdot, t)=P_tf$ be a positive solution to the heat
equation $\partial_t u=Lu$ with $u(\cdot, 0)=f$, $f$ is a positive and measurable function on $M$. Let
\begin{eqnarray*}
H_{K}(f, t)=D_K(t)\int_M (P_t(f\log f)-P_tf\log P_tf )d\mu,
\end{eqnarray*}
where $D_0(t)={1\over t}$ and $D_{K}(t)={1\over |1-e^{-2Kt}|}$ for $K\neq 0$.Then, for all $K\in \mathbb{R}$,
\begin{eqnarray*}
{d\over dt}H_{K}(f, t)\leq 0,\ \ \ \forall t>0,
\end{eqnarray*}
and for all $K\in \mathbb{R}$ and $t>0$, we have
\begin{eqnarray*}
{d^2\over dt^2}H_K(t)+2K\coth(Kt) {d\over dt}H_K(t)
\leq - 2D_K(t)\int_M |\nabla^2\log P_tf|^2P_tfd\mu. \label{KK2}
\end{eqnarray*}
Define the $W$-entropy by the revised Boltzmann entropy formula
\begin{eqnarray*}
W_K(f, t)=H_K(f, t)+{\sinh(2Kt)\over 2K}{d\over dt}H_K(f, t).
\end{eqnarray*}
Then, for all $K\in \mathbb{R}$, and for all $t>0$, we have
\begin{eqnarray*}
{d\over dt}W_{K}(f, t)&=&-{\sinh(2Kt)\over K}D_K(t)\int_M  |\nabla^2
\log P_tf|^2 P_tf d\mu\nonumber\\
& & -{\sinh(2Kt)\over K}D_K(t)\int_M (Ric(L)-K)(\nabla\log P_tf, \nabla\log P_tf)P_tfd\mu.\label{WWWW}
\end{eqnarray*}
In particular, for all $K\in \mathbb{R}$, we have
\begin{eqnarray*}
{d\over dt}W_{K}(f, t)\leq 0,\ \ \ \forall t>0.
\end{eqnarray*}
\end{theorem}

\begin{theorem} \label{Th-W2} Let $M$ be a complete Riemannian manifold with a fixed metric and potential $(g, \phi)$. Suppose that $(M, g)$ satisfies the bounded geometry condition and $\phi\in C^4(M)$ with $\nabla\phi\in C^3_b(M)$. Let $u$ be the heat kernel of the Witten Laplacian $L=\Delta-\nabla\phi\cdot\nabla$. Let
\begin{eqnarray*}
H_{m, K}(u, t)=-\int_M u\log u d\mu-\Phi_{m, K}(t),
\end{eqnarray*}
where $\Phi_{m, K}\in C((0, \infty), \mathbb{R})$ satisfies
\begin{eqnarray*}
\Phi_{m, K}'(t)={m\over 2t}e^{4Kt} ,\ \ \ \forall t>0.
\end{eqnarray*}
Define the $W$-entropy by the Boltzmann formula
\begin{eqnarray*}
W_{m, K}(u, t)={d\over dt}(tH_{m, K}(u, t)).
\end{eqnarray*}
Then
\begin{eqnarray*}
{d\over dt}W_{m, K}(u, t)
&=&-2t\int_M \left[\left|\nabla^2\log u+\left({K\over 2}+{1\over 2t}\right)g\right|^2+(Ric_{m, n}(L)+Kg)(\nabla\log u, \nabla \log u)\right] ud\mu\\
& &\hskip1.5cm -{2t\over m-n}\int_M \left|\nabla \phi\cdot \nabla\log  u-{(m-n)(1+Kt)\over 2t}\right|^2ud\mu\\
& &\hskip3cm -{m\over 2t}\left[e^{4Kt}(1+4Kt)-(1+Kt)^2\right].
\end{eqnarray*}
In particular, if $Ric_{m, n}(L)\geq -K$, then, for all $t\geq 0$,
we have
\begin{eqnarray*}
{d\over dt}W_{m, K}(u, t)\leq -{m\over
2t}\left[e^{4Kt}(1+4Kt)-(1+Kt)^2\right].
\end{eqnarray*}
Moreover, the equality holds at some time $t=t_0>0$ if and only if
$M$ is a quasi-Einstein manifold, i.e., $Ric_{m, n}(L)=-Kg$, and the
potential function $f=-\log u$ satisfies the shrinking soliton
equation with respect to $Ric_{m, n}(L)$, i.e.,
\begin{eqnarray*}
Ric_{m, n}(L)+2\nabla^2f={g\over t},
\end{eqnarray*}
and moreover
\begin{eqnarray*}
\nabla \phi\cdot \nabla f=-{(m-n)(1+Kt)\over 2t}.
\end{eqnarray*}
\end{theorem}

We would like to point out that, Theorem \ref{WCDK} and  Theorem
\ref{Th-W2} are new even in the case $\phi$ is a constant, $m=n$ and
$L=\Delta$ is the usual Laplace-Beltrami operator on complete
Riemannian manifolds with Ricci curvature bounded from below by a
negative constant. In this case, Theorem \ref{Th-W2} can be formulated as follows.

\begin{theorem} \label{Th-W2b} Let $(M, g)$ be a complete Riemannian manifold with bounded geometry
condition. Then
\begin{eqnarray*}
{d\over dt}W_{n, K}(u, t)
&=&-2t\int_M u\left[\left|\nabla^2\log u+\left({K\over 2}+{1\over 2t}\right)g\right|^2+(Ric+Kg)(\nabla\log u, \nabla \log u)\right] ud\mu\\
& &\hskip3cm -{n\over 2t}\left[e^{4Kt}(1+4Kt)-(1+Kt)^2\right].
\end{eqnarray*}
In particular, if $Ric\geq -K$, then, for all $t\geq 0$, we have
\begin{eqnarray*}
{d\over dt}W_{n, K}(u, t)\leq -{n\over
2t}\left[e^{4Kt}(1+4Kt)-(1+Kt)^2\right].
\end{eqnarray*}
Moreover, the equality holds at some time $t=t_0>0$ if and only if
$M$ is an Einstein manifold, i.e., $Ric=-Kg$, and the potential
function $f=-\log u$ satisfies the shrinking soliton equation, i.e.,
\begin{eqnarray*}
Ric+2\nabla^2f={g\over t}.
\end{eqnarray*}
\end{theorem}

In Subsection $3.3$, we will extend Theorem \ref{WCDK} and Theorem
\ref{Th-W2} to time dependent Witten Laplacian on compact Riemannian manifolds with a $K$-super Perelman  Ricci flow. For details, see Theorem \ref{WCDK3}  and
Theorem \ref{Th-W3}.

In Section $4$,  we prove  the Li-Yau Harnack inequality and the Li-Yau-Hamilton Harnack inequality for positive solutions to the heat equation $\partial_t u=Lu$ of the time dependent Witten Laplacian on compact Riemannian manifolds equipped with variants 
of the $(K, m)$-super Ricci flow. 

\begin{theorem}\label{LYHSRF-A}  Let $(M, g(t), \phi(t), t\in [0, T])$ be a compact Riemannian manifold with a family of time dependent metrics $g(t)$ and potentials $\phi(t)\in C^2(M)$, $t\in [0, T]$. Let $u$ be a positive solution to the heat equation $\partial_t u=Lu$. Let $\partial_t g=2h$ and $\alpha>1$. Suppose that  $(M, g(t), \phi(t), t\in [0, T])$  satisfies the backward $(\alpha, K, m)$-super Ricci flow 
\begin{eqnarray}
{1\over 2}(1-\alpha)\partial_t g+Ric_{m, n}(L)\geq -Kg,\label{mmm1}
\end{eqnarray}
and assume that $A^2=\max\limits \left[ |h|^2+{({\rm Tr}h)^2\over m-n}\right]<\infty$ and $B=\max\limits |S|<\infty$, where
\begin{eqnarray*}
S(\cdot)=2h(\nabla \phi, \cdot\rangle-\langle 2{\rm div} h-\nabla {\rm Tr}_g h+\nabla \phi_t, \cdot\rangle+{2 {\rm Tr} h\over m-n}\langle \nabla\phi, \cdot\rangle.
\end{eqnarray*}
Then for any $\gamma>0$ and for all $t\in (0, T]$, we have

\begin{eqnarray*}
{|\nabla u|^2\over u^2}-\alpha {\partial_t u\over u}\leq {m\alpha^2\over 4t}\left[1+\sqrt{1+{T^2\over m}\left(4A^2+{(2K+\gamma)^2\over (\alpha-1)^2}+{2B^2\over \gamma}\right)}\right].
\end{eqnarray*}
In the case $B=0$,  for all $t\in (0, T]$, we have
\begin{eqnarray*}
{|\nabla u|^2\over u^2}-\alpha {\partial_t u\over u}\leq {m\alpha^2 \over 4t}\left[1+\sqrt{1+{T^2\over m}\left(4A^2+{4K^2\over (\alpha-1)^2}\right)}\right].
\end{eqnarray*}
In particular, in the case $A=B=0$ and $Ric_{m, n}(L)\geq 0$, we have the Li-Yau Harnack inequality
\begin{eqnarray*}
{|\nabla u|^2\over u^2}-{\partial_t u\over u}\leq {m\over 2t}.
\end{eqnarray*}
\end{theorem}

\begin{theorem}\label{LYHHSRF-A}  Let $(M, g(t), t\in [0, T])$ be a compact Riemannian manifold with a family of time dependent metrics $g(t)$ and potentials $\phi(t)\in C^2(M)$, $t\in [0, T]$. Let $u$ be a positive solution to the heat equation $\partial_t u=Lu$. Let $\partial_t g=2h$. 
 Suppose that, for all $t\in (0, T]$,
\begin{eqnarray}
e^{-4Kt} \left(h+Ric_{m, n}(L)+Kg\right)-e^{-2Kt}h\geq  \alpha_K(t)g,
\end{eqnarray}
where $\alpha_K(t)$ is a real valued function in time $t$, and
\begin{eqnarray*}
A^2=\max\limits \left[ |h|^2+{({\rm Tr}h)^2\over m-n}\right]<\infty,\ \ \ \ \ \ B=\max\limits |S|<\infty,
\end{eqnarray*}
where
\begin{eqnarray*}
S(\cdot)=\left\langle {2{\rm Tr}h\over m-n} \nabla \phi-2 {\rm div}h-\nabla {\rm Tr}_g h+\nabla \partial_t \phi, \cdot\right\rangle +2h(\nabla\phi, \cdot).
\end{eqnarray*}
Then, for any $\gamma>0$ and $t\in [0, T]$, we have
\begin{eqnarray*}
{|\nabla u|^2\over u^2}-e^{2Kt}{\partial_t u\over u}  \leq  {me^{4Kt}\over 2t} \left[1+\sqrt{{A^2T^2\over m}+\max\limits_{t\in [0, T]}{t^2(2\alpha_K(t)-\gamma)^2\over 4e^{-4Kt}(1-e^{-2Kt})^2}+\max\limits_{t\in [0, T]}{t^2e^{-4Kt}B^2\over 2m\gamma} } \right].
\end{eqnarray*}
In the case $B=0$, we have
\begin{eqnarray*}
{|\nabla u|^2\over u^2}-e^{2Kt}{\partial_t u\over u}  \leq  {me^{4Kt}\over 2t} \left[1+\sqrt{{A^2T^2\over m}+\max\limits_{t\in [0, T]}{t^2\alpha^2_K(t)\over e^{-4Kt}(1-e^{-2Kt})^2} } \right].
\end{eqnarray*}
and if $\alpha_K(t)=0$, i.e., if
\begin{eqnarray}
e^{-4Kt} (h+Ric_{m, n}(L)+K)-e^{-2Kt}h\geq  0,
\end{eqnarray}
we have
\begin{eqnarray*}
{|\nabla u|^2\over u^2}-e^{2Kt}{\partial_t u\over u}  \leq  {me^{4Kt}\over 2t} \left[1+{TA\over \sqrt{m}} \right].
\end{eqnarray*}
In particular, when $A=B=0$, and $Ric_{m, n}(L)\geq -K$,  we recapture Hamilton's Harnack inequality \cite{H1}
\begin{eqnarray*}
{|\nabla u|^2\over u^2}-e^{2Kt}{\partial_t u\over u}  \leq {me^{4Kt}\over 2t}.
\end{eqnarray*}
\end{theorem}

Theorem \ref{LYHSRF-A} and Theorem \ref{LYHHSRF-A} can be also extended to complete Riemannian manifolds equipped with variants of the $(K, m)$-super Ricci flow.  For details, see Section $5$. 

The rest of this paper is organized as follows.  In Section $2$, we prove Theorem \ref{Thm-RLSI}, Theorem \ref{Thm0}, Corollary \ref{cor1}, Theorem \ref{HLYH} and Corollary \ref{cor2}.  In Section $3$,  we prove Theorem \ref{WCDK}, Theorem \ref{Th-W2} and extend Theorem \ref{WCDK} and Theorem \ref{Th-W2} to time dependent Witten Laplacian on compact Riemannian manifolds with a $K$-super Perelman  Ricci flow, see Theorem \ref{WCDK3} and Theorem \ref{Th-W3}. In Section $4$, we prove Theorem \ref{LYHSRF-A} and Theorem \ref{LYHHSRF-A}. In Section $5$,  we extend Theorem \ref{LYHSRF-A}  and Theorem \ref{LYHHSRF-A} to complete Riemannian manifolds equipped with variants of the $(K, m)$-super Ricci flows.  

\section{Log-Sobolev inequalities and Harnack inequalities for Witten Laplacian}

\subsection{Log-Sobolev inequalities on $K$-super Perelman  Ricci flow}

In this subsection, we modify the semigroup argument due to Bakry and Ledoux \cite{BL} to prove  the equivalence between the $K$-super Perelmam Ricci flow equation and two logarithmic Sobolev inequalities for  the time dependent Witten Laplacian.

{\it Proof}.  Let $P_{s, t}$ be the heat semigroup of the time dependent weighted
Laplacian on $L^2(M, \mu)$, i.e., for any $s\in [0, T]$, $u(t, \cdot):=P_{s, t}f(\cdot)$ is the unique solution of
the heat equation $\partial_t u=Lu$ in $L^2(M, \mu)$ on $[s, T]$ with $u(s, \cdot)=f$. Let 
\begin{eqnarray*}
h(s, t)=e^{2Kt}P_{s+T-t, T}\left({|\nabla P_{s, s+T-t}f|^2\over P_{s, s+T-t}f}\right),
\ \ \ t\in [s, T].
\end{eqnarray*}
Note that, at time $T-t+s$, the generalized Bochner formula implies
\begin{eqnarray}
\label{Boh1}
(\partial_t+L){|\nabla u|^2\over u}={2\over u}|\nabla^2
u-u^{-1}\nabla u\otimes \nabla u|^2+2u^{-1}\left({1\over 2}{\partial
g\over
\partial t}+Ric(L)\right)(\nabla u, \nabla u).
\end{eqnarray}
Hence
\begin{eqnarray*}
\partial_t h(s, t)&=&2K h(s, t)+e^{2Kt}P_{s+T-t, T}\left[\left({\partial \over
\partial t}+L\right)\left({|\nabla P_{s, s+T-t}f|^2\over
P_{s, s+T-t}f}\right)\right]\\
&=&2K h(s, t)+e^{2Kt} P_{s+T-t, T}\left[{2\over u}|\nabla^2 u-u^{-1}\nabla u\otimes
\nabla u|^2+2u^{-1}\left({1\over 2}{\partial g\over
\partial t}+Ric(L)\right)(\nabla u, \nabla u) \right]\\
 &\geq& 2e^{2Kt} P_{s+T-t, T}\left[u^{-1}\left({1\over 2}{\partial g\over
\partial t}+Ric(L)+K\right)(\nabla u, \nabla u)\right].
\end{eqnarray*}
If $(g(t), \phi(t))$ is a $(K, \infty)$-super Ricci flow, i.e., if  ${1\over 2}{\partial g\over
\partial t}+Ric(L)+K \geq 0$, then
$$
\partial_t h(s, t) \geq 0.
$$
Thus, $t\rightarrow h(s, t)$ is increasing on $[s, T]$. This yields,
for all $t\in (s, T)$,
\begin{eqnarray*}\label{f1}
e^{2Ks}{|\nabla P_{s, T}f|^2\over P_{s, T}f}\leq e^{2Kt}P_{s+T-t, T}\left({|\nabla
P_{s, s+T-t}f|^2\over P_{s, s+T-t}f}\right)\leq e^{2KT}P_{s,T}\left({|\nabla
f|^2\over f}\right).
\end{eqnarray*}
Notice that
\begin{eqnarray*}
{d\over dt}P_{s+T-t, T}(P_{s, s+T-t}f\log
P_{s, s+T-t}f)&=&P_{s+T-t, T}\left((L_{s+T-t}+\partial_t)(P_{s, s+T-t}f\log
P_{s, s+T-t}f)\right)\\
&=&P_{s+T-t, T}\left({|\nabla P_{s, s+T-t} f|^2\over P_{s, s+T-t}f}\right).
\end{eqnarray*}
Therefore,
\begin{eqnarray*}
P_{s,T}(f\log f)-P_{s,T}f\log P_{s, T} f&=&\int_s^T {d\over dt}P_{s+T-t, T}(P_{s, s+T-t}f\log
P_{s, s+T-t}f)dt\\
&=&\int_s^T P_{s+T-t, T}\left({|\nabla P_{s, s+T-t} f|^2\over P_{s, s+T-t}f}\right)dt\\
&\leq& \frac{1}{2K}(e^{2K(T-s)} - 1)P_{s, T}\left({|\nabla f|^2\over f}\right).
\end{eqnarray*}
Thus the logarithmic Sobolev inequality holds on complete Riemannian manifolds equipped with a $K$-super Perelman  Ricci flow
\begin{eqnarray}\label{LSI2}
P_{s,T}(f\log f)-P_{s,T}f\log P_{s,T} f\leq {e^{2K(T - s)}-1\over 2K}P_{s,T}\left({|\nabla f|^2\over f}\right).
\end{eqnarray}
Similarly to the above proof of the logarithmic Sobolev inequality $(\ref{LSI2})$, we have
\begin{eqnarray*}
P_{s, T}(f\log f)-P_{s, T}f\log P_{s,T} f &=& \int_s^T {d\over dt}P_{s+T-t, T}(P_{s, s+T-t}f\log P_{s, s+T-t}f)dt\\
&=&\int_s^T P_{s+T-t, T}\left({|\nabla P_{s, s+T-t} f|^2\over P_{s, s+T-t}f}\right)dt\\
&\geq&\int_s^T e^{2K(s-t)}{|\nabla P_{s,T} f|^2\over
P_{s,T} f}dt\\
&=&{1-e^{2K(s-T)}\over 2K}{|\nabla P_{s,T} f|^2\over P_{s, T} f}.
\end{eqnarray*}
Thus,  for all $T>0$, $f\in C_b(M)$ with $f>0$, the reversal  logarithmic Sobolev inequality holds
\begin{eqnarray}\label{rLSI2}
{|\nabla P_{s,T} f|^2\over P_{s, T} f}\leq {2K\over 1-e^{2K(s-T)} }\left(P_{s,T}(f\log
f)-P_{s,T}f\log P_{s, T}f\right).
\end{eqnarray}
Changing $T$ by $t$, we see that $\eqref{LSI2}$ and $\eqref{rLSI2}$
hold for $P_{s, t}$ for all $0 \leq s< t \leq T$.

On the other hand, if for all $0 \leq s < t \leq T$, the log-Sobolev
inequality holds for $P_{s, t}$, then applying $(\ref{LSI2})$ to
$1+\varepsilon f$ and letting $\varepsilon\rightarrow 0$, we can
obtain the Poincar\'e inequality
\begin{eqnarray*}
P_{s,t}f^{2}- (P_{s,t}f)^{2} \leq \frac{1}{K}(e^{2K(t-s)} - 1)P_{s, t}\left({|\nabla f|^2}\right).
\end{eqnarray*}
Set
\begin{eqnarray*}
w(s, t)=P_{s,t}f^{2}- (P_{s,t}f)^{2}-\frac{1}{K}(e^{2K(t-s)} - 1)P_{s, t}\left({|\nabla f|^2}\right).
\end{eqnarray*}
Then $w(s, t)\leq 0$ for all $0\leq s<t<T$. Notice that when $s = t$, we have $w(s, s)=0$. Hence, for all $s<t$, $\partial_s w(s, t)\leq 0$. Now
\begin{eqnarray*}
\partial_{s} w(s, t) &=& -P_{s, t}L_sf^{2} + 2P_{s, t}f P_{s, t}L_sf + 2e^{2K(t-s)}P_{s, t} (|\nabla f|^{2}) \\
& &\ \ \ - \frac{e^{2K(t-s)} - 1} {K} [P_{s, t}(-L(|\nabla f|^{2})+\partial_s|\nabla f|^{2}_{g(s)})].
\end{eqnarray*}
Thus,  at $s =t$, we have $\left. \partial_s w(s,
t)\right|_{s=t}=0$. On the other hand, as for all $s<t$, $\partial_s
w(s, t)\leq 0$, we can derive that $\partial^2_s w(s, t)\leq 0$ for
all $0\leq s<t<T$. By calculation, at $s=t$, we have
\begin{eqnarray*}
\partial_{s}^2 w(s, t) &=&-\partial_s(P_{s, t}L_sf^{2}) + 2\partial_s(P_{s, t}fP_{s, t}L_sf) - 4Ke^{2K(t-s)}P_{s, t}(|\nabla f|^{2}) \\
& &\ \ + 2e^{2K(t-s)}\partial_sP_{s, t}(|\nabla f|^{2}) + 2e^{2K(t-s)}P_{s, t}(-L_s(|\nabla f|^{2}) + \partial_s|\nabla f|^{2}_{g(s)})\\
&=& L^{2}_sf^{2}  - \partial_s(L_sf^{2}) - 2(L_s f)^{2} - 2fL^{2}_{s}f + 2f\partial_sL_sf\\\
& &\ \ \ \ - 4K(|\nabla f|^{2}) - 4L_s(|\nabla f|^{2})
+4\partial_s|\nabla f|_{g(s)}^2.
\end{eqnarray*}
Note that
\begin{eqnarray*}
L_s^2 f^2&=&L_s(2fL_sf+2|\nabla f|_{g(s)}^2)\\
&=&2L_s fL_s f+2fL_s^2f+4\nabla f\cdot\nabla L_sf+2L_s|\nabla
f|_{g(s)}^2,\\
\partial_s L_s f^2&=&\partial_s (2fL_s f+2|\nabla f|_{g(s)}^2)=2f\partial_s L_s f+2\partial_s |\nabla f|_{g(s)}^2,
\end{eqnarray*}
and
\begin{eqnarray*}
\partial_s |\nabla f|_{g(s)}^2=-\partial_s g(s)(\nabla f, \nabla f),
\end{eqnarray*}
from which and the Bochner formula, we see that at $s=t$,
\begin{eqnarray*}
\partial_{s}^2 w(s, t)
&=&-2L_s|\nabla f|_{g(s)}^2+4\nabla f\cdot \nabla L_s f-4K|\nabla
f|_{g(s)}^2+2\partial_s|\nabla f|_{g(s)}^2\\
&=&-4\left.\left[\|{\rm Hess} f|_{g(s)}^2+\left({1\over 2}\partial_s
g(s)+Ric(L_s)+K\right)(\nabla f, \nabla f)\right]\right|_{s=t}.
\end{eqnarray*}
Taking $f$ to be normal coordinate functions near $x$ on $(M,
g(t))$, we derive that, at any time $t \in [0, T]$,
$$
{1\over 2}{\partial g\over
\partial t} + Ric(L_t) + K \geq 0.
$$
This completes the proof of Theorem \ref{Thm-RLSI}.

\begin{remark}{\rm
Indeed, we can further prove
the Poincar\'e inequality, the reversal Poinca\'e
inequality as well as Bakry-Ledoux's Gromov-L\'evy isoperimetric
inequality on the super Ricci flow $(\ref{RFK})$. To save the length of the paper, we will do these in a forthcoming paper. In \cite{St1, St2,
St3}, Sturm introduced the super Ricci flow on metric measure
space, and proved the equivalence between the super Ricci flow and
the Poincar\'e inequality.
}
\end{remark}

\subsection{Hamilton's Harnack inequality for Witten Laplacian with $CD(K, \infty)$ condition}

In this subsection, we prove Hamilton's Harnack inequality for the time dependent Witten Laplacian on complete Riemannian manifolds with a $K$-super Perelman  Ricci flow.

{\it Proof of Theorem \ref{Thm0}}. We modify  the method used in \cite{Li13}. Let $t\in [0, T)$  and $s\in [0, T-t]$. Using the reversal logarithmic Sobolev inequality and the fact $0<f\leq A$, we have
\begin{eqnarray*}
{|\nabla P_{s, s+t} f|^2\over P_{s, s+t} f}&\leq& {2K\over
1-e^{-2Kt}}\left(P_{s, s+t}(f\log f)-P_{s, s+t}f \log P_{s, s+t}f\right)\\
&\leq& {2K\over
1-e^{-2Kt}}\left(P_{s, s+t}(f\log A)-P_{s, s+t}f \log P_{s, s+t}f\right).
\end{eqnarray*}
Thus
\begin{eqnarray*}
|\nabla \log P_{s, s+t} f|^2\leq {2K\over 1-e^{-2Kt}}\log (A/P_{s, s+t}f).
\end{eqnarray*}
Using ${1\over 1-e^{-x}}\leq 1+{1\over x}$ for $x\geq 0$, we have
\begin{eqnarray*}
|\nabla \log P_{s, s+t} f|^2 \leq \left(2K+{1\over t}\right)\log
(A/P_{s, s+t}f).
\end{eqnarray*}
In particular, for $s=0$, we have 
\begin{eqnarray*}
{|\nabla u|^2\over u^2} \leq \left(2K+{1\over t}\right)\log
(A/u).
\end{eqnarray*}
The proof of Theorem \ref{Thm0} is completed. \hfill $\square$

\medskip

{\it Proof of Corollary \ref{cor1}}. Let $l(x, t)=\log{A/u(x, t)}$. Then the differential Harnack  inequality $(\ref{HH})$ in Theorem \ref{Thm0} implies
\begin{eqnarray*}
|\nabla\sqrt{l(x, t)}|={1\over 2}{|\nabla l(x, t)|\over \sqrt{l(x, t)}}\leq {1\over 2}\sqrt{2K\over 1-e^{-2Kt}}.
\end{eqnarray*}
Fix $x, y\in M$ and integrate along a geodesic linking $x$ and $y$, the above inequality yields
\begin{eqnarray*}
\sqrt{\log{A/u(x, t)}}\leq \sqrt{\log{A/u(y, t)}}+{1\over 2}\sqrt{2K\over 1-e^{-2Kt}}d(x, y).
\end{eqnarray*}
Combining this with the elementary inequality
\begin{eqnarray*}
(a+b)^2\leq (1+\delta)a^2+(1+\delta^{-1})b^2,
\end{eqnarray*}
we can derive the desired Harnack inequality for $u$. \hfill $\square$

\subsection{The LYH Harnack inequality for Witten Laplacian with $CD(K, m)$ condition}

In this subsection we prove Theorem \label{HLYH}, i.e,   the Li-Yau-Hamilton Harnack inequality for the positive solution to the heat equation associated with the Witten Laplacian on complete Riemannian manifold with fixed metric and potential. Indeed, by the generalized Bochner-Weitzenb\"ock formula, we have
\begin{eqnarray}
(L-\partial_t){|\nabla u|^2\over u}={2\over u}\left|\nabla^2
u-{\nabla u\otimes \nabla u\over u}\right|^2+{2\over u}Ric(L)(\nabla
u, \nabla u).\label{GB2}
\end{eqnarray}
Taking trace in the first quantity on the right hand side, we can
derive
\begin{eqnarray}
(L-\partial_t){|\nabla u|^2\over u}\geq {2\over nu}\left|\Delta
u-{|\nabla u|^2\over u}\right|^2+{2\over u}Ric(L)(\nabla u, \nabla
u).\label{GB3}
\end{eqnarray}
Applying the inequality
\begin{eqnarray*}
(a+b)^2\geq {a^2\over 1+\alpha}-{b^2\over \alpha}
\end{eqnarray*}
to $a=\partial_t u-{|\nabla u|^2\over u}$, $b=\nabla \phi\cdot\nabla u$, and $\alpha={m-n\over n}$, we have
\begin{eqnarray}
(L-\partial_t){|\nabla u|^2\over u}\geq {2\over m u}\left|\partial_t
u-{|\nabla u|^2\over u}\right|^2+{2\over u} Ric_{m, n}(L)(\nabla u,
\nabla u). \label{GB4}
\end{eqnarray}
Hence, under the condition $Ric_{m, n}(L)\geq -K$, it holds
\begin{eqnarray}
(L-\partial_t){|\nabla u|^2\over u} \geq  {2\over m
u}\left|\partial_t u-{|\nabla u|^2\over u}\right|^2-{2K|\nabla
u|^2\over u}.\label{GB5}
\end{eqnarray}
Let
\begin{eqnarray*}
h={\partial u \over \partial t}-e^{-2Kt}{|\nabla u|^2\over u}+e^{2Kt}{m\over 2t}u.
\end{eqnarray*}
Then $\lim\limits_{t\rightarrow 0^+}h(t)=+\infty$, and
\begin{eqnarray*}
(\partial_t-L)h\geq {2\over m u}e^{-2Kt}\left|\partial_t u-{|\nabla u|^2\over u}\right|^2-e^{2Kt}{m\over 2t^2}u.
\end{eqnarray*}

We now prove that $h\geq 0$ on $M\times \mathbb{R}^+$. In compact case, suppose that $h$ attends its minimum at some $(x_0, t_0)$ and $h(x_0, t_0)<0$. Then, at $(x_0, t_0)$, it holds
\begin{eqnarray*}
{\partial h\over \partial t}\leq 0,\ \ \Delta h\geq 0, \ \nabla h=0.
\end{eqnarray*}
Thus at $(x_0, t_0)$, $(\partial_t-L)h\leq 0$.  On the other hand, at this point, we have
\begin{eqnarray*}
0\leq e^{2Kt}{m\over 2t}u<e^{-2Kt}{|\nabla u|^2\over u}-{\partial u\over \partial t}\leq {|\nabla u|^2\over u}-{\partial u\over \partial t},
\end{eqnarray*}
and hence
\begin{eqnarray*}
(\partial_t-L) h>0.
\end{eqnarray*}
This finishes the proof of Theorem \ref{HLYH} in compact case.

In complete non-compact case, let $f=\log u$, and let
\begin{eqnarray*}
F=te^{-2Kt}(e^{-2Kt}|\nabla f|^2-f_t)=te^{-4Kt}|\nabla f|^2-te^{-2Kt}f_t.
\end{eqnarray*}
Obviously, $F(0, x)\equiv 0$. We shall prove that
\begin{eqnarray*}
F\leq {m\over 2}.
\end{eqnarray*}
By direct calculation
\begin{eqnarray*}
LF&=&te^{-4Kt}L|\nabla f|^2-te^{-2Kt}Lf_t\\
\partial_t F&=&\partial_t (te^{-4Kt}|\nabla f|^2-te^{-2Kt}f_t)\\
&=&(1-4Kt)e^{-4Kt}|\nabla f|^2+(2Kt-1)e^{-2Kt}f_t+te^{-4Kt}\partial_t |\nabla f|^2-te^{-2Kt}f_{tt},
\end{eqnarray*}
we have
\begin{eqnarray*}
(L-\partial_t)F&=&te^{-4Kt}(L-\partial_t)|\nabla f|^2-te^{-2Kt}(L-\partial_t)f_t\\
& &\ \ \ +(4Kt-1)e^{-4Kt}|\nabla f|^2-(2Kt-1)e^{-2Kt}f_t.
\end{eqnarray*}
By the generalized Bochner formula, it holds
\begin{eqnarray*}
(L-\partial_t)|\nabla f|^2=2|\nabla^2 f|^2+2Ric(L)(\nabla f, \nabla f)-4\nabla^2f(\nabla f, \nabla f).
\end{eqnarray*}
Note that
\begin{eqnarray*}
Lf_t&=&L\left({Lu\over u}\right)
={L^2u\over u}-2\langle\nabla Lu, {\nabla u\over u^2}\rangle+Lu\left(-{Lu\over u^2}+2{|\nabla u|^2\over u^3}\right),\\
\partial_t f_t&=&\partial_t\left({Lu\over u}\right)={L^2u\over u}-{|Lu|^2\over u^2},
\end{eqnarray*}
which yields
\begin{eqnarray*}
(L-\partial_t)f_t&=& 2{Lu|\nabla u|^2\over u^3} -2\langle\nabla Lu, {\nabla u\over u^2}\rangle\\
&=&-4\nabla^2 f(\nabla f, \nabla f)-2\langle\nabla Lf, \nabla f\rangle.
\end{eqnarray*}
Hence
\begin{eqnarray*}
(L-\partial_t)F
&=& 2te^{-4Kt}[|\nabla^2 f|^2 + 2(e^{2Kt}-1)\nabla^2 f(\nabla f, \nabla f)]\\
& &+2te^{-4Kt}Ric(L)(\nabla f, \nabla f) +2 te^{-2Kt}\langle\nabla Lf, \nabla f\rangle\\
& &\ \ \ +(4Kt-1)e^{-4Kt}|\nabla f|^2-(2Kt-1)e^{-2Kt}(Lf+|\nabla f|^2).
\end{eqnarray*}
Now
\begin{eqnarray*}
F &=& te^{-4Kt}|\nabla f|^2-te^{-2Kt}f_t= te^{-4Kt}(1 - e^{2Kt})|\nabla f|^2 - te^{-2Kt}Lf,\\
\langle \nabla F, \nabla f \rangle &=& 2te^{-4Kt}(1 - e^{2Kt})\nabla^2f(\nabla f, \nabla f) - te^{-2Kt}\langle\nabla Lf, \nabla f\rangle.
\end{eqnarray*}
Therefore
\begin{eqnarray*}
(L-\partial_t)F
&=&  2te^{-4Kt}|\nabla^2 f|^2 -2\langle \nabla F, \nabla f \rangle\\
& &\ \ +2te^{-4Kt}\left(Ric(L)(\nabla f, \nabla f)+K|\nabla
f|^2\right)+ \frac{(2Kt-1)}{t}F.
\end{eqnarray*}
By \cite{Li05}, we have
$$
|\nabla^2 f|^2 \geq \frac{1}{n}|\Delta f|^2 \geq \frac{1}{m}|L f|^2
-{1\over m-n}\nabla\phi\otimes\nabla \phi(\nabla f, \nabla f).
$$
Thus
\begin{eqnarray*}
(L-\partial_t)F&\geq& 2te^{-4Kt}{|Lf|^2\over m}-2\langle \nabla F, \nabla f \rangle\\
& &\ \ +2te^{-4Kt}\left(Ric_{m, n}(L)(\nabla f, \nabla f)+K|\nabla
f|^2\right)+ \frac{(2Kt-1)}{t}F\\
&\geq& \frac{2te^{-4Kt}}{m}\left[\frac{(te^{-2Kt}(e^{-2Kt}
-1)|\nabla f|^2 - F)^2}{t^2e^{-4Kt}}\right]-2\langle \nabla F, \nabla f \rangle +  \frac{(2Kt-1)}{t}F\\
&\geq& \frac{2[te^{-2Kt}(e^{-2Kt} -1)|\nabla f|^2 - F]^2}{mt}-
2\langle \nabla F, \nabla f \rangle +  \frac{(2Kt-1)}{t}F.
\end{eqnarray*}

Similarly to \cite{Li05}, let $\eta$ be a $C^2$-function on $[0,―infty)$ such that $\eta=1$ on $[0, 1]$ and $\eta=0$ on $[2,―infty)$, with $-C_1\eta^{1/2}(r)\leq \eta'(r)\leq 0$, and $\eta^{''}(r)\geq C_2$, where $C_, C_2>0$ are two constants. Let $\rho(x)=d(o, x)$ and define $\psi(x)=\eta({\rho(x)/R})$.
Since $\rho$ is Lipschitz on the complement of the cut locus of $o$, $\psi$ is a Lipschitz function with support
in $B(o, 2R)\times [0,―infty)$. As explained in Li and Yau [43], an argument of Calabi allows us
to apply the maximum principle to $\psi F$. Let $(x_0, t_0)\in M\times [0, T]$ be
a point where $\psi F$ achieves the maximum. Then, at $(x_0, t_0)$,
\begin{eqnarray*}
\partial_t(\psi F)\geq 0, \ \Delta(\psi F)\leq 0, \ \nabla(\psi F)=0.
\end{eqnarray*}
This yields
\begin{eqnarray*}
(L-\partial_t)(\psi F)=\Delta (\psi F)-\nabla\phi\cdot\nabla(\psi
F)-\partial_t(\psi F)\leq 0.
\end{eqnarray*}
Similarly to \cite{Li05}, we have
\begin{eqnarray*}
(L-\partial_t)(\psi F)&=&\psi(L-\partial_t)F+(L\psi)F+2\nabla\psi \cdot\nabla F\\
&\geq&\psi(L-\partial_t)F-A(R)F+2\nabla\psi\cdot\nabla F\\
&\geq &\psi(L-\partial_t)F -A(R)F+2\langle\nabla \psi, \nabla(\psi
F)\rangle\psi^{-1}-2F|\nabla\psi|^2\psi^{-1}.
\end{eqnarray*}where we use
\begin{eqnarray*}
L\psi\geq -A(R):=-{C_1\over
R}(m-1)\sqrt{K}\coth(\sqrt{K}R)-{C_2\over R^2},
\end{eqnarray*}
and for some constant $C_3>0$
$${|\nabla \psi|^2\over \psi}\leq {C_3\over R^2}.$$
Let $C(n, K, R)={C_1\over
R}(m-1)\sqrt{K}\coth(\sqrt{K}R)+{C_2+C_3\over R^2}$. At the point $(x_0, t_0)$, we have
\begin{eqnarray*}
0 &\geq &\psi(L-\partial_t)F -(A(R) + 2|\nabla\psi|^2\psi^{-1})F\\
&\geq & \psi\left[\frac{2[te^{-2Kt}(e^{-2Kt} -1)|\nabla f|^2 - F]^2}{mt}- 2\langle \nabla F, \nabla f \rangle +  \frac{(2Kt-1)}{t}F\right] -C(n, K, R)F\\
&\geq& \psi\frac{2}{mt}F^2 + \psi\frac{4e^{-2Kt}(1 - e^{-2Kt})|\nabla f|^2}{m}F + 2F\langle \nabla \psi, \nabla f \rangle +  \left[(2K -\frac{1}{t})\psi -C(n, K, R)\right]F\\
&\geq& \psi\frac{2}{mt}F^2 + \psi\frac{4e^{-2Kt}(1 - e^{-2Kt})|\nabla f|^2}{m}F - 2F|\nabla \psi||\nabla f| +  \left[(2K -\frac{1}{t})\psi -C(n, K, R)\right]F\\
&\geq& \psi\frac{2}{mt}F^2 + \psi\frac{4e^{-2Kt}(1 - e^{-2Kt})|\nabla f|^2}{m}F - 2{C_2 \over R}F\psi^{1/2}|\nabla f| +  \left[(2K -\frac{1}{t})\psi -C(n, K, R)\right]F.
\end{eqnarray*}
Multiplying by $t$ on both sides, and using the Cauchy-Schwartz inequality, we get
\begin{eqnarray*}
0 
&=& \psi\frac{2}{m}F^2 + tF\left[\psi\frac{4e^{-2Kt}(1 - e^{-2Kt})|\nabla f|^2}{m} - 2{C_2 \over R}\psi^{1/2}|\nabla f|\right] +  [(2Kt - 1)\psi - C(n, K, R)t]F\\
&\geq& \psi\frac{2}{m}F^2 +  \left[(2Kt - 1)\psi - C(n, K, R)t -
{C_2 m\over 4e^{-2Kt}(1 - e^{-2Kt})R^2}t\right]F.
\end{eqnarray*}
Notice that the above calculation is done at the point $(x_0, t_0)$.
Since $\psi F$ reaches its maximum at this point, we can assume that
$\psi F(x_0, t_0) > 0$. Thus
\begin{eqnarray*}
0 
\geq \frac{2}{m}(\psi F)^2  -\left[1 + C(n, K, R)t + {C_2 m\over
4e^{-2Kt}(1 - e^{-2Kt})R^2}t\right](\psi F),
\end{eqnarray*}
which yields that, for any $(x, t) \in B_{R} \times [0, T]$,
\begin{eqnarray*}
F(x,t) &\leq& (\psi F)(x_0, t_0) \leq \frac{m}{2}\left[1 + C(n, K, R)t_0 + {C_2m \over 4e^{-2Kt_0}(1 - e^{-2Kt_0})R^2}t_0\right]\\
&\leq& \frac{m}{2} \left[1 + C(n, K, R)T + \max\limits_{t\in [0,
T]}{C_2 mt\over 4e^{-2Kt}(1 - e^{-2Kt})R^2}\right].
\end{eqnarray*}
Let $R \rightarrow \infty$, we obtain
$$
F \leq \frac{m}{2}.
$$
\hfill $\square$

\medskip

{\it Proof of Corollary \ref{cor2}}.  The proof is as the same as the one of Corollary 2.2 in \cite{H1}. For the completeness we reproduce it as follows. Let $l(x, t)=\log u(x, t)$. Then the Li-Yau-Hamilton Harnack inequality is equivalent to
\begin{eqnarray}
{\partial l\over \partial t}-e^{-2Kt}|\nabla l|^2+e^{2Kt}{m\over 2t}\geq 0.\label{HHHH}
\end{eqnarray}
Let $\gamma: [0, T]\rightarrow M$ be a geodesic with reparametrization by arc length $s: [\tau, T]\rightarrow [0, T]$ so that $\gamma(s(\tau))=x$ and $\gamma(s(T))=y$. Let $S(t)={d\gamma(s(t))\over dt}=\dot\gamma(s(t))\dot s(t)$. Then
$|\dot\gamma(s(t))|=1$. Integrating along $\gamma(s(t))$ from $t=\tau$ to $t=T$, we have
\begin{eqnarray*}
l(y, T)-l(x, \tau)=\int_\tau^T \left[{\partial l\over \partial t}+\nabla l\cdot S\right]dt.
\end{eqnarray*}
By the Cauchy-Schwartz inequality
\begin{eqnarray*}
e^{-2Kt}|\nabla l|^2+{1\over 4}e^{2Kt}|S|^2\geq \nabla l\cdot S
\end{eqnarray*}
From this and $(\ref{HHHH})$ we obtain
\begin{eqnarray*}
l(y, T)-l(x, \tau)\geq -{1\over 4}\int_\tau^Te^{2Kt}|S|^2dt-\int_\tau^T {m\over 2t}e^{2Kt}dt.
\end{eqnarray*}
Note that $d(x, y)=\int_\tau^T |S|dt=\int_\tau^T ds(t)$.
Choosing $s(t)=a[e^{-2K\tau}-e^{-2Kt}]$, with
$$a={d(x, y)\over e^{-2K\tau}-e^{-2KT}},$$ we have
\begin{eqnarray*}
l(y, T)-l(x, \tau)&\geq &-{1\over 4}\int_t^T e^{2Kt}{\dot s}^2(t)dt-\int_\tau^T {m\over 2t}e^{2Kt}dt\\
&=&-{Kd^2(x, y)\over 2(e^{-2K\tau}-e^{-2KT})}-\int_\tau^T {m\over 2t}e^{2Kt}dt.
\end{eqnarray*}
Note that $\int_\tau^T {e^{2Kt}\over t}dt\leq \log\left({T\over \tau}\right)+e^{2KT}-e^{2K\tau}$. Thus
\begin{eqnarray*}
\log u(y, T)-\log u(x, \tau)\geq -{Kd^2(x, y)\over 2(e^{-2K\tau}-e^{-2KT})}-{m\over 2}\left[\log\left({T\over \tau}\right)+e^{2KT}-e^{2K\tau}\right].
\end{eqnarray*}
Using ${1\over 1-e^{-x}}\leq {1+x\over x}$, we can derive the desired estimate. \hfill $\square$

\subsection{Hamilton's second order  estimates for Witten Laplacian with $CD(K, m)$ condition}

In \cite{H1}, Hamilton also proved that, on compact Riemannian manifolds with $Ric\geq -K$,  there exists a constant $C$ depending only on $n$ and $K$, such that for any positive solution of the heat equation $\partial_t u=\Delta u$ with $0<u\leq A$ and $t\in [0, 1]$, it holds
\begin{eqnarray}
t\Delta u\leq C u\left[1+\log({A/u})\right].\label{HHHH1}
\end{eqnarray}
Indeed, Hamilton \cite{H1} proved the following estimate
\begin{eqnarray}
{\Delta u\over u} +{|\nabla u|^2\over u^2}\leq {K\over 1-e^{-Kt}}\left[n+\log(A/u)\right],\ \ \ \forall\ t\geq 0.\label{HHHH2}
\end{eqnarray}

In this subsection, we extend Hamilton's second order estimate $(\ref{HHHH2})$ to positive solutions of the heat equation $\partial_t u=Lu$ for the Witten Laplacian on complete Riemannian manifolds with $CD(K, m)$-condition. More precisely,  we prove the following

\begin{theorem}\label{HH3} Let $m\geq n$ and $K\geq 0$ be two constants. Let $M$ be a complete Riemannian manifold with a $C^2$-potential such that the $CD(-K, m)$ condition holds
$$Ric_{m, n}(L)\geq -K.$$ Let $u$ be a
positive solution of the heat equation
\begin{eqnarray*}
\partial_t u=Lu
\end{eqnarray*}
with $A=\sup\limits\{u(x, t), (x, t)\in M\times [0, T]\}<\infty$. Then
\begin{eqnarray}
{L u\over u}+{|\nabla u|^2\over u^2}\leq {K\over 1-e^{-Kt}}\left[m+4\log(A/u)\right],\ \ \ \forall\ t\in [0, T].\label{HHHH3}
\end{eqnarray}
In particular, for $t\in [0, T]$, we have
\begin{eqnarray}
{Lu\over u}\leq \left(K+{1\over t}\right)\left[m+\log(A/u)\right].\label{HHHH4}
\end{eqnarray}
\end{theorem}
{\it Proof}. Let $\psi(t)={1-e^{-Kt}\over K}$. Then $\psi'+K\psi=1$. Let $h=\psi \left[Lu+{|\nabla u|^2\over u}\right]-u[m+4\log(A/u)]$.
By $(\ref{GB4})$ , under the assumption $Ric_{m, n}(L)\geq -K$ we have
\begin{eqnarray*}
(\partial_t-L){|\nabla u|^2\over u}\leq -{2\over m u}\left|L u-{|\nabla u|^2\over u}\right|^2+2K{|\nabla u|^2\over u},
\end{eqnarray*}
which yields
\begin{eqnarray*}
(\partial_t-L)h\leq -{2\psi\over mu}\left|L u-{|\nabla u|^2\over u}\right|^2+\psi' \left[L u-{|\nabla u|^2\over u}\right]-2{|\nabla u|^2\over u}.
\end{eqnarray*}
By analogue of Hamilton\cite{H1}, we can verify that
\begin{eqnarray*}
{\partial h\over \partial t}\leq Lh \ \ \ {\rm whenever}\ \ h\geq 0.
\end{eqnarray*}
Indeed, we can verify this by examining three cases:\\

$(i)$ If $L u\leq {|\nabla u|^2\over u}$, then $(\partial_t-\Delta )h\leq 0$ since $\psi'\geq 0$.

$(ii)$ If ${|\nabla u|^2\over u}\leq Lu\leq 3{|\nabla u|^2\over u}$, then $(\partial_t-L)h\leq 0$ since $\psi'\leq 1$.

$(iii)$ If $3{|\nabla u|^2\over u}\leq Lu$, then whenever $h\geq 0$, we have
\begin{eqnarray*}
2\left[L u-{|\nabla u|^2\over u}\right]\geq L u+{|\nabla u|^2\over u}={h\over \psi}+{mu+4u\log(A/u)\over \psi} \geq {mu\over \psi},
\end{eqnarray*}
which yields, since $\psi'\leq 1$, we have
\begin{eqnarray*}
(\partial_t-L)h\leq (\psi'-1)\left[Lu-{|\nabla u|^2\over u}\right]-2{|\nabla u|^2\over u}\leq 0.
\end{eqnarray*}

Note that $h\leq 0$ at $t=0$. By the maximum principle, we conclude that $h\leq 0$ for all $t\geq 0$. Thus
\begin{eqnarray*}
{Lu\over u}+{|\nabla u|^2\over u^2}\leq {K\over 1-e^{-Kt}}\left[m+4\log(A/u)\right].
\end{eqnarray*}
This completes the proof of Theorem \ref{HH3}. \hfill $\square$

\section{The $W$-entropy formula for Witten Laplacian}

Recall that, Perelman \cite{P1} introduced the notion of the
$W$-entropy and proved its monotonicity along the conjugate heat
equation associated to the Ricci flow. In \cite{N1, N2}, Ni proved
the monotonicity of the $W$-entropy for the heat equation of the
usual Laplace-Beltrami operator on complete Riemannian manifolds
with non-negative Ricci curvature. In \cite{Li12, Li13},  the second
author of this paper proved the $W$-entropy formula  and its monotonicity and rigidity theorems for the heat
equation of the Witten Laplacian on complete Riemannian manifolds satisfying the $CD(0, m)$ condition
and gave a probabilistic interpretation of the $W$-entropy for the
Ricci flow. In \cite{LL13}, we gave a new proof of the $W$-entropy
formula obtained in \cite{Li12} for the Witten Laplacian  by
using Ni's $W$-entropy formula $(\ref{entropy-2})$ to the
Laplace-Beltrami operator on $M\times S^{m-n}$ equipped with a suitable warped product Riemannian metric,
and further proved the monotonicity of the $W$-entropy for the heat
equation of the time dependent Witten Laplacian on compact
Riemannian manifolds equipped with the super Ricci flow with respect
to the $m$-dimensional Bakry-Emery Ricci curvature.

During the past years, many people have asked the following very
natural problem to us.

\begin{problem} How to define the $W$-entropy  for the heat equation associated with the Witten
Laplacian on complete Riemannian manifolds satisfying the $CD(K, m)$ condition for $K\in \mathbb{R}$ and $m\in [n, \infty]$? Can we establish the monotonicity and rigidity  theorems for the $W$-entropy associated with the Witten
Laplacian on complete Riemannian manifolds satisfying  general curvature-dimension condition? What happens on manifolds with time dependent metrics and potentials?
\end{problem}

In this section, we give the answer to this fundamental problem.

\subsection{$W$-entropy for Witten Laplacian with $CD(K, \infty)$ condition}

In this subsection, based on the reversal  logarithmic Sobolev inequality  on complete Riemannian manifolds with fixed metrics and potentials, which is due to Bakry and Ledoux \cite{BL},  we introduce the  $W$-entropy and  prove
 the $W$-entropy formula for the Witten Laplacian on complete Riemannian manifolds with fixed metric and potential satisfying the $CD(K, \infty)$ condition.

Let $C_0(t)={1\over t}$, and for $K\neq 0$, $C_K(t)={2K\over e^{2Kt}-1}$.  Let  $D_0(t)={1\over t}$, $D_{K}(t)={2|K|\over |1-e^{-2Kt}|}$.
Then $D_K'(t)=-C_K(t)D_K(t)$ for all $K\in \mathbb{R}$ and $t>0$.
We first introduce the revised Boltzmann-Shannon entropy
\begin{eqnarray*}
H_{K}(f, t)=D_K(t)\int_M (P_t(f\log f)-P_t f\log P_tf)d\mu,
\end{eqnarray*}
where $f$ is a positive and measurable function on $M$. Based on the gradient estimates of the positive solution to the heat
equation on complete manifolds with bounded geometry condition (see
\cite{Li12, Li13}), by direct calculation and using the integration
by parts formula , we can prove
\begin{eqnarray}
{d\over d t} H_{K}(f, t)&=&C_K(t) D_K(t)\int_M  \left(P_tf\log P_tf-P_t(f\log f)\right)d\mu+D_K(t)\int_M {|\nabla P_t f|^2\over P_t f} d\mu\nonumber\\
&=&D_K(t)\int_M \left[{|\nabla P_tf|^2\over P_tf}+C_K(t)(P_tf\log P_tf-P_t(f\log f))\right]d\mu.\label{HHHH1}
\end{eqnarray}
Under the condition  $Ric(L)\geq K$, by the reversal  logarithmic Sobolev inequality due to Bakry and Ledoux \cite{BL}, for all $t>0$, we have
\begin{eqnarray}
{|\nabla P_tf|^2\over P_tf}\leqq C_K(t)(P_t(f\log f)-P_tf\log P_tf). \label{BLRLSI}
\end{eqnarray}
Hence, for all $K\in \mathbb{R}$, we have
\begin{eqnarray*}
{d\over d t} H_{K}(f, t)\leq 0,\ \ \ \forall\ t>0.
\end{eqnarray*}
Taking the time derivative on the both sides of $(\ref{HHHH1})$, we have
\begin{eqnarray*}
{d^2\over dt^2}H_{K}(f, t)&=&-C_K(t)D_K(t)\left[\int_M {|\nabla P_tf|^2\over P_tf}+C_K(t)P_tf \log P_tf-K(t)P_t(f\log f)\right]d\mu\\
& & \hskip0.5cm +D_K(t)\left[{d\over dt}\int_M {|\nabla
P_tf|^2\over P_tf}d\mu-C_K(t)\int_M {|\nabla P_tf|^2\over P_tf}d\mu\right]\\
& &\hskip1cm +D_K(t){d\over dt}C_K(t)\int_M (P_tf\log P_tf-P_t(f\log
f))d\mu.
\end{eqnarray*}
By Bakry and Emery \cite{BE} and Li \cite{Li12, Li13}, we have
\begin{eqnarray*}
{d\over dt}\int_M {|\nabla P_tf|^2\over P_tf}d\mu=-2\int_M |\nabla^2\log
P_tf|^2P_tfd\mu-2\int_M Ric(L)(\nabla \log P_tf, \nabla \log P_tf)P_tfd\mu.
\end{eqnarray*}
Thus
\begin{eqnarray}
{d^2\over dt^2}H_{K}(f, t)
&=&-C_K(t){d\over dt}H_K(t)-2D_K(t) \int_M  |\nabla^2
\log P_tf|^2 P_tf d\mu\nonumber\\
& &\hskip0.5cm -D_K(t)\int_M (2Ric(L)+C_K(t))(\nabla \log P_tf, \nabla\log P_tf)P_tfd\mu\nonumber\\
& &\hskip1cm +D_K(t){d\over dt}C_K(t)\int_M (P_tf\log P_tf-P_t(f\log
f))d\mu\label{KH1}.
\end{eqnarray}
Note that, for all $K\in \mathbb{R}$, under the condition $Ric(L)\geq K$, we have
\begin{eqnarray*}
2Ric(L)+C_K(t)\geq 2K+{2K\over e^{2Kt}-1}={2K\over 1-e^{-2Kt}},
\end{eqnarray*}
and
\begin{eqnarray*}
{d\over dt}C_K(t)={d\over dt}{2K\over e^{2Kt}-1}=-{2K\over 1-e^{-2Kt}}C_k(t).
\end{eqnarray*}
Substituting these into $(\ref{KH1})$,  a simple calculation yields, for all $K\in \mathbb{R}$, and for all $t>0$,
\begin{eqnarray*}
{d^2\over dt^2}H_{K}(f, t)\leq -2K\coth(Kt){d\over dt}H_K(t)-2D_K(t) \int_M  |\nabla^2
\log P_tf|^2 P_tf d\mu.\label{HHH}
\end{eqnarray*}
Indeed, from $(\ref{KH1})$, we can prove
\begin{eqnarray*}
{d^2\over dt^2}H_{K}(f, t)&=&-2K\coth(Kt){d\over dt}H_K(t)-2D_K(t) \int_M  |\nabla^2
\log P_tf|^2 P_tf d\mu\nonumber\\
& &\hskip1cm -2D_K(t)\int_M (Ric(L)-K)(\nabla \log P_tf, \nabla \log P_tf)P_tfd\mu.\label{HHH2}
\end{eqnarray*}

Let $\alpha_K: (0, \infty)\rightarrow (0, \infty)$ be a $C^1$-smooth function. Define the $W$-entropy by the revised Boltzmann entropy formula
\begin{eqnarray*}
W_{K}(f, t):={1\over \dot \alpha_K(t)}{d\over dt}(\alpha_K(t)H_K(f, t))= H_K+{\alpha_K\over \dot \alpha_k} \dot H_K.
\end{eqnarray*}
Set $\beta_K={\alpha_K\over \dot\alpha_K}$. Then
\begin{equation*}
{d\over dt}W_{K}(f, t)=\beta_K(\ddot{H}_K+\frac{1+\dot{\beta}_K}{\beta_K}H_K).
\end{equation*}
Solving the ODE
\begin{equation*}
\frac{1+\dot{\beta_K}}{\beta_K}=2K\coth(Kt),
\end{equation*}
we can take \begin{eqnarray*}
\beta_K(t)={\sinh(2Kt)\over 2K},
\end{eqnarray*}
and hence
\begin{eqnarray*}
\alpha_K(t)=K\tanh(Kt).
\end{eqnarray*}
This yields
\begin{eqnarray*}
W_K(f, t)=H_K(f, t)+{\sinh(2Kt)\over 2K}{d\over dt}H_K(f, t),
\end{eqnarray*}
and
\begin{eqnarray*}
{d\over dt}W_{K}(f, t)&=&-{\sinh(2Kt)\over K}D_K(t)\int_M  |\nabla^2
\log P_tf|^2 P_tf d\mu\nonumber\\
& & -{\sinh(2Kt)\over K}D_K(t)\int_M (Ric(L)-K)(\nabla\log P_tf, \nabla\log P_tf)P_tfd\mu.\label{WWWW}
\end{eqnarray*}
In particular, when $Ric(L)\geq  K$, then for all $t>0$, we have
\begin{eqnarray*}
{d\over dt}W_{K}(f, t)
\leq 0.
\end{eqnarray*}
This finishes the proof of Theorem \ref{WCDK}. \hfill $\square$

\subsection{$W$-entropy for Witten Laplacian with $CD(K, m)$ condition}

In this subsection, based on Theorem \ref{HLYH},  we introduce the $W$-entropy and prove
 the $W$-entropy formula and a rigidity theorem, i.e.,  Theorem \ref{Th-W2}, for the Witten Laplacian on complete Riemannian manifolds with fixed metrics and potentials satisfying the $CD(K, m)$ condition for $K\in \mathbb{R}$ and $m\in [n, \infty)$.

Let $M$ be a complete Riemannian manifold with bounded geometry
condition, $\phi\in C^3(M)$ be such that $\nabla\phi\in C_b^2(M)$. Following \cite{P1, N1, Li12, Li11, LL13},
we define
\begin{eqnarray*}
H_{m, K}(u, t)=-\int_M u\log u d\mu-\Phi_{m, K}(t),
\end{eqnarray*}
where $\Phi_{m, K}\in C((0, \infty), \mathbb{R})$ satisfies
\begin{eqnarray*}
\Phi_{m, K}'(t)={m\over 2t}e^{4Kt}, \ \ \ \forall t>0.
\end{eqnarray*}

\begin{proposition}
Let $M$ be a complete Riemannian manifold with bounded geometry
condition, $\phi\in C^3(M)$ be such that $\nabla\phi\in C_b^2(M)$.
Then, under the condition $Ric_{m, n}(L)\geq -K$, we have
\begin{eqnarray*}
{d\over dt}H_{K, m}(u, t)\leq 0.
\end{eqnarray*}
\end{proposition}
{\it Proof}. By the entropy dissipation formula (see \cite{Li12,
Li13}, and using $\int_M
\partial_t u d\mu=\int_M Lu d\mu=0$, we have
\begin{eqnarray}
{d\over dt}H_{m, K}(u, t)&=&\int_M \left[{|\nabla u|^2\over u^2}-{m\over 2t}e^{4Kt}\right] ud\mu\label{HmK1}\\
&=&\int_M \left[{|\nabla u|^2\over u^2}-{m\over 2t}e^{4Kt}-e^{2Kt}{\partial_t u\over u}\right] ud\mu.\label{HmK2}
\end{eqnarray}
By the Li-Yau-Hamilton Harnack inequality in Theorem \ref{HLYH}, we have
\begin{eqnarray*}
{d\over dt}H_{m, K}(u, t)\leq 0.
\end{eqnarray*}
\hfill $\square$

\begin{proposition}\label{prop2} Under the same condition as in
Theorem \ref{Th-W2}, we have
\begin{eqnarray*}
{d^2\over dt^2}H_{m, K}(u, t)=-2\int_M [|\nabla^2\log u|^2+Ric(L)(\nabla\log u, \nabla \log u)]ud\mu-\left({2mK\over t}-{m\over 2t^2}\right)e^{4Kt}.
\end{eqnarray*}
\end{proposition}
{\it Proof}. Based on the gradient estimates of the positive
solution to the heat equation on complete manifolds with bounded
geometry condition (see \cite{Li12, Li13}), we have
\begin{eqnarray*}
{d\over dt}\int_M {|\nabla u|^2\over u}d\mu=-2\int_M [|\nabla^2\log u|^2+Ric(L)(\nabla\log u, \nabla \log u)]ud\mu.
\end{eqnarray*}
Combining this with $(\ref{HmK1})$, Proposition \ref{prop2} follows.
\hfill $\square$

\medskip

Following Perelman \cite{P1}, Ni \cite{N1} and \cite{Li12, Li11, Li13, LL13}, we
introduce the $W$-entropy for the heat equation $(\ref{Heat1})$ of the Witten Laplacian as
follows
\begin{eqnarray*}
W_{m, K}(u, t)={d\over dt}( tH_{m, K}(u, t)).
\end{eqnarray*}
By the entropy dissipation formula, we have
\begin{eqnarray*}
W_{m, K}(u, t)&=&\int_M \left[t(|\nabla \log u|^2-\Phi'_{m,
K}(t))-\log u-\Phi_{m, K}(t)\right]ud\mu\\
&=&\int_M \left[t(2L(-\log u)-|\nabla \log u|^2)-\log u-\Phi_{m,
K}(t)-\Phi_{m, K}'(t)\right]ud\mu.
\end{eqnarray*}

We are now in a position to state the main result of this section, i.e., Theorem \ref{Th-W2}.

\begin{theorem} \label{Th-W2b} Let $M$ be a complete Riemannian manifold, $\phi\in C^2(M)$. Suppose that $(M, g)$ satisfies the bounded geometry condition and $\phi\in C^4(M)$ with $\nabla\phi\in C^3_b(M)$.
Then
\begin{eqnarray*}
{d\over dt}W_{m, K}(u, t)
&=&-2t\int_M \left[\left|\nabla^2\log u+\left({K\over 2}+{1\over 2t}\right)g\right|^2+(Ric_{m, n}(L)+Kg)(\nabla\log u, \nabla \log u)\right] ud\mu\\
& &\hskip1.5cm -{2t\over m-n}\int_M \left|\nabla \phi\cdot \nabla\log  u-{(m-n)(1+Kt)\over 2t}\right|^2ud\mu\\
& &\hskip3cm -{m\over 2t}\left[e^{4Kt}(1+4Kt)-(1+Kt)^2\right].
\end{eqnarray*}
In particular, if $Ric_{m, n}(L)\geq -K$, then, for all $t\geq 0$,
we have
\begin{eqnarray*}
{d\over dt}W_{m, K}(u, t)\leq -{m\over
2t}\left[e^{4Kt}(1+4Kt)-(1+Kt)^2\right].
\end{eqnarray*}
Moreover, the equality holds at some time $t=t_0>0$ if and only if
$M$ is a quasi-Einstein manifold, i.e., $Ric_{m, n}(L)=-Kg$, and the
potential function $f=-\log u$ satisfies the shrinking soliton
equation with respect to $Ric_{m, n}(L)$, i.e.,
\begin{eqnarray*}
Ric_{m, n}(L)+2\nabla^2f={g\over t},
\end{eqnarray*}
and moreover
\begin{eqnarray*}
\nabla \phi\cdot \nabla f=-{(m-n)(1+Kt)\over 2t}.
\end{eqnarray*}
\end{theorem}
{\it Proof}. By $(\ref{HmK1})$ and Proposition \ref{prop2}, we have
\begin{eqnarray*}
{d\over dt}W_{m, K}(u, t)
&=&-2t\left[\int_M [|\nabla^2\log u|^2+Ric(L)(\nabla\log u, \nabla \log u)]ud\mu+\left({mK\over t}-{m\over 4t^2}\right)e^{4Kt}\right]\\
& &\ \ +2\int_M \left[{|\nabla u|^2\over u^2}-{m\over 2t}e^{4Kt}\right] ud\mu.
\end{eqnarray*}
Note that
\begin{eqnarray*}
\left|\nabla^2\log u+\left({e^{2Kt}\over
2t}+a(t)\right)g\right|^2=|\nabla^2\log u|^2+2\left({e^{2Kt}\over
2t}+a(t)\right)\Delta \log u+n\left({e^{2Kt}\over 2t}+a(t)\right)^2.
\end{eqnarray*}
By direct calculation, we have
\begin{eqnarray*}
{d\over dt}W_{m, K}(u, t)
&=&-2t\int_M \left|\nabla^2\log u+\left({e^{2Kt}\over 2t}+a(t)\right)g\right|^2 u d\mu\\
& &-2t\int_M \left(Ric_{m, n}(L)+\left(2a(t)-{1-e^{2Kt}\over t}\right)g\right)(\nabla\log u, \nabla \log u) ud\mu\\
& &\ \ \ +2nt\left({e^{2Kt}\over 2t}+a(t)\right)^2-{me^{4Kt}\over 2t}-2mK e^{4Kt}\\
& &\ \ \ +2(e^{2Kt}+2ta(t))\int_M \nabla \phi\cdot \nabla\log  u \
ud\mu-2t\int_M {|\nabla \phi\cdot\nabla \log u|^2 \over m-n}ud\mu.
\end{eqnarray*}
Let $a(t)$ be chosen such that $2a(t)-{1-e^{2Kt}\over t}=K$.
Then
\begin{eqnarray*}
{d\over dt}W_{m, K}(u, t)
&=&-2t\int_M \left[\left|\nabla^2\log u+\left({K\over 2}+{1\over 2t}\right)g\right|^2+(Ric_{m, n}(L)+Kg)(\nabla\log u, \nabla \log u)\right] ud\mu\\
& &\ \ \ +2nt\left({1\over 2t}+{K\over 2}\right)^2-{me^{4Kt}\over 2t}-2mK e^{4Kt}\\
& &\ \ \ +2(1+Kt)\int_M \nabla \phi\cdot \nabla\log  u \ ud\mu-2t\int_M {|\nabla \phi\cdot\nabla \log u|^2 \over m-n}ud\mu.
\end{eqnarray*}
Combining this with
\begin{eqnarray*}
& &{1\over m-n}\int_M \left|\nabla\phi\cdot\nabla\log u-{(m-n)(1+Kt)\over 2t}\right|^2ud\mu\\
&=&{(m-n)(1+Kt)^2\over 4t^2}-{1+Kt\over t}\int_M \nabla \phi\cdot \nabla\log  u \ ud\mu+\int_M {|\nabla \phi\cdot\nabla \log u|^2 \over m-n}ud\mu,
\end{eqnarray*}
and noting that
\begin{eqnarray*}
& & 2nt\left({1\over 2t}+{K\over 2}\right)^2-{me^{4Kt}\over 2t}-2mK e^{4Kt}+{(m-n)(1+Kt)^2\over 2t}\\
& &\hskip3cm ={m\over 2t}\left[(1+Kt)^2-e^{4Kt}(1+4Kt)\right],
\end{eqnarray*}
we can derive the desired  $W$-entropy formula. The rest of the
proof is obvious.  \hfill $\square$

In particular, taking $m=n$, $\phi\equiv 0$ and $g$ is a fixed Riemannian metric, we obtain the following $W$-entropy formula for the heat equation of the Laplace-Beltrami operator on Riemannian manifolds, which extends Ni's result in \cite{N1} for $K=0$.

\begin{theorem} \label{Th-W2c} Let $(M, g)$ be a complete Riemannian manifold with bounded geometry
condition. Let $u$ be the fundamental solution to the heat equation
$\partial_t u=\Delta u$. Then
\begin{eqnarray*}
{d\over dt}W_{n, K}(u, t)
&=&-2t\int_M \left[\left|\nabla^2\log u+\left({K\over 2}+{1\over 2t}\right)g\right|^2+(Ric+Kg)(\nabla\log u, \nabla \log u)\right] ud\mu\\
& &\hskip3cm -{n\over 2t}\left[e^{4Kt}(1+4Kt)-(1+Kt)^2\right].
\end{eqnarray*}
In particular, if $Ric\geq -K$, then, for all $t\geq 0$, we have
\begin{eqnarray*}
{d\over dt}W_{n, K}(u, t)\leq -{n\over
2t}\left[e^{4Kt}(1+4Kt)-(1+Kt)^2\right].
\end{eqnarray*}
Moreover, the equality holds at some time $t=t_0>0$ if and only if
$M$ is an Einstein manifold, i.e., $Ric=-Kg$, and the potential
function $f=-\log u$ satisfies the shrinking soliton equation, i.e.,
\begin{eqnarray*}
Ric+2\nabla^2f={g\over t}.
\end{eqnarray*}
\end{theorem}

To end this subsection, let us remark that in our previous paper \cite{LL13} we introduced another $W$-entropy functional for  the heat equation associated with the Witten Laplacian satisfying the $CD(K, m)$-condition as follows

\begin{eqnarray*}
\widetilde{W}_{m, K}(u)={d\over dt}(t\widetilde{H}_{m, K}(u)),
\end{eqnarray*}
where
\begin{eqnarray*}
\widetilde{H}_{m, K}(u)=-\int_M u\log ud\mu-{m\over 2t}(1+\log(4\pi t))-{mKt\over 2}(1+{1\over 6}Kt),
\end{eqnarray*}
and we proved that
\begin{eqnarray*}
{d\over dt}\widetilde{W}_{m, K}(u)&=&-2t\int_M \left[\left|\nabla^2\log u+\left({K\over 2}+{1\over 2t}\right)g\right|^2+(Ric_{m, n}(L)+Kg)(\nabla\log u, \nabla \log u)\right] ud\mu\\
& &\hskip1.5cm -{2t\over m-n}\int_M \left|\nabla \phi\cdot \nabla\log  u-{(m-n)(1+Kt)\over 2t}\right|^2ud\mu.
\end{eqnarray*}
Indeed, letting $\Psi_{m, K}(t)=\Phi_{m, K}(t)-{m\over 2t}(1+\log(4\pi t))-{mKt\over 2}(1+{1\over 6}Kt)$, we have
\begin{eqnarray*}
\widetilde{W}_{m, K}(u)-W_{m, K}(u)={d\over dt}(t\Psi_{m, K}(t)),
\end{eqnarray*}
and
\begin{eqnarray*}
{d\over dt}(\widetilde{W}_{m, K}(u)-W_{m, K}(u))={d^2\over dt^2}(t\Psi_{m, K}(t))={m\over 2t}\left[e^{4Kt}(1+4Kt)-(1+Kt)^2\right].
\end{eqnarray*}

\subsection{$W$-entropy for Witten Laplacian on $K$-super Perelman  Ricci flow}

In this subsection, we extend the $W$-entropy formula to the time
dependent Witten Laplacian on compact Riemannian manifolds with $K$-super Perelman  Ricci flow.

Let $(M, g(t), \phi(t), t\in [0, T])$ be a complete Riemannian manifold with a family of time dependent metrics $g(t)$ and potentials $\phi(t)$. Let
$$L=\Delta_{g(t)}-\nabla_{g(t)}\phi(t)\cdot\nabla_{g(t)}$$
be the time dependent Witten Laplacian on $(M, g(t), \phi(t))$. Let
$$d\mu(t)=e^{-\phi(t)}dvol_{g(t)}.$$
Suppose that
\begin{eqnarray}
{\partial \phi\over \partial t}={1\over 2}{\rm Tr} {\partial g\over
\partial t}.\label{PPPPP}
\end{eqnarray}
Then $\mu(t)$ is indeed independent of $t\in [0, T]$, i.e.,
\begin{eqnarray*}
{\partial \mu(t)\over \partial t}=0, \ \ \ t\in [0, T].
\end{eqnarray*}

We now state the main results of this subsection, which extend
Theorem \ref{WCDK} and Theorem \ref{Th-W2} to the time dependent
Witten Laplacian on compact Riemannian manifolds with $K$-super Perelman  Ricci flow.

\begin{theorem}\label{WCDK3} Let $(M, g(t), t\in [0, T])$ be a compact Riemannian manifold with a family of
metrics $g(t)$, and $\phi\in C^{2, 1}(M\times [0, T])$. Suppose that
$(\ref{PPPPP})$ holds and
$${1\over 2}{\partial g\over \partial t}+Ric(L)\geq K,$$ where $K\in \mathbb{R}$ is
a constant. Let $u(\cdot, t)=P_tf$ be a positive solution to the heat
equation $\partial_t u=Lu$ with $u(\cdot, 0)=f$, $f$ is a positive and measurable function on $M$. Define
\begin{eqnarray*}
H_{K}(f, t)=D_K(t)\int_M (P_t(f\log f)-P_tf\log P_tf )d\mu,
\end{eqnarray*}
where $D_0(t)={1\over t}$ and $D_{K}(t)={1\over |1-e^{-2Kt}|}$ for $K\neq 0$.Then, for all $K\in \mathbb{R}$,
\begin{eqnarray*}
{d\over dt}H_{K}(f, t)\leq 0,\ \ \ \forall t\in (0, T],
\end{eqnarray*}
and for all $K\in \mathbb{R}$ and $t\in (0, T]$, we have
\begin{eqnarray*}
{d^2\over dt^2}H_K(t)+2K\coth(Kt) {d\over dt}H_K(t)
\leq - 2D_K(t)\int_M |\nabla^2\log P_tf|^2P_tfd\mu.
\end{eqnarray*}
Define the $W$-entropy by the revised Boltzmann entropy formula
\begin{eqnarray*}
W_K(f, t)=H_K(f, t)+{\sinh(2Kt)\over 2K}{d\over dt}H_K(f, t).
\end{eqnarray*}
Then, for all $K\in \mathbb{R}$, and for all $t\in (0, T]$, we have
\begin{eqnarray*}
{d\over dt}W_{K}(f, t)&=&-{\sinh(2Kt)\over K}D_K(t)\int_M  |\nabla^2
\log P_tf|^2 P_tf d\mu\nonumber\\
& & -{\sinh(2Kt)\over K}D_K(t)\int_M \left({1\over 2}{\partial g\over \partial t}+Ric(L)-K\right)(\nabla\log P_tf, \nabla\log P_tf)P_tfd\mu.\label{WWWW2}
\end{eqnarray*}
In particular, for all $K\in \mathbb{R}$, we have
\begin{eqnarray*}
{d\over dt}W_{K}(f, t)\leq 0,\ \ \ \forall t\in (0, T].
\end{eqnarray*}
\end{theorem}
{\it Proof}. By Li-Li \cite{LL13}, we have
\begin{eqnarray*}
{d\over dt}\int_M {|\nabla P_tf|^2\over P_tf}d\mu=-2\int_M |\nabla^2\log
P_tf|^2P_tfd\mu-2\int_M \left({1\over 2}{\partial g\over \partial
t}+Ric(L)\right)(\nabla \log P_tf, \nabla \log P_tf)P_tfd\mu.
\end{eqnarray*} The rest of the proof is similar the one of Theorem \ref{WCDK}.
\hfill $\square$

\begin{theorem} \label{Th-W3} Let $(M, g(t), \phi(t), t\in [0, T])$ be a
compact Riemannian manifold with a family of time dependent metrics
and potentials $(g(t), \phi(t), t\in [0, T])$ satisfying
$(\ref{PPPPP})$. Let $m\geq n$ and $K\geq 0$ be two constants which
are independent of $t\in [0, T]$. Let $u$ be a positive solution to
the heat equation ${\partial_t u}=Lu$. Let
\begin{eqnarray*}
W_{m, K}(u, t)={d\over dt}( tH_{m, K}(u, t)),
\end{eqnarray*}
where
\begin{eqnarray*}
H_{m, K}(u, t)=-\int_M u\log u d\mu-\Phi_{m, K}(t).
\end{eqnarray*}
Then, for all $t\in [0, T]$, we have
\begin{eqnarray*}
{d\over dt}W_{m, K}(u, t)
&=&-2t\int_M \left|\nabla^2\log u+\left({K\over 2}+{1\over 2t}\right)g\right|^2ud\mu-{m\over 2t}\left[e^{4Kt}(1+4Kt)-(1+Kt)^2\right]\\
& &\hskip0.5cm -2t\int_M \left({1\over 2}{\partial g\over \partial t}+Ric_{m, n}(L)+Kg\right)(\nabla\log u, \nabla \log u)ud\mu\\
& &\hskip1.5cm -{2t\over m-n}\int_M \left|\nabla \phi\cdot
\nabla\log  u-{(m-n)(1+Kt)\over 2t}\right|^2ud\mu.
\end{eqnarray*}
In particular, if $(M, g(t), \phi(t), t\in [0, T])$ is a $K$-super
Perelman Ricci flow with respect to the $m$-dimensional
Bakry-Emery Ricci curvature $Ric_{m, n}(L)$,
$${1\over 2}{\partial g\over \partial t}+Ric_{m, n}(L)\geq -Kg,$$
we have
\begin{eqnarray*}
{d\over dt}W_{m, K}(u, t)\leq -{m\over
2t}\left[e^{4Kt}(1+4Kt)-(1+Kt)^2\right].
\end{eqnarray*}
Moreover, the equality holds at some time $t=t_0>0$ if and only if
$(M, g(t), t\in [0, T])$ is a quasi-Ricci flow, i.e.,
\begin{eqnarray*}
{1\over 2}{\partial g\over \partial t}&=&-Ric_{m, n}(L)-Kg,\\
{\partial\phi\over \partial t}&=&-R-\Delta \phi+{|\nabla\phi|^2\over
m-n}-nK,
\end{eqnarray*}
and the potential function $f=-\log u$ satisfies
\begin{eqnarray*}
2\nabla^2f &=&\left({1\over t}+K\right)g,\\
\nabla \phi\cdot \nabla f&=&-{(m-n)(1+Kt)\over 2t}.
\end{eqnarray*}
\end{theorem}
{\it Proof}. The proof is similar to the one of Theorem \ref{Th-W2}. See \cite{LL13} for the case $K=0$. \hfill $\square$

\medskip

Similarly to the end of Section $3.2$, we can reformulate Theorem \ref{Th-W3} in terms of $\widetilde {W}_{m, K}$. See \cite{LL13}.  

\section{The Li-Yau and the Li-Yau-Hamilton Harnack inequalities on compact super-Ricci flows}

In this section we prove the Li-Yau Harnack inequality and the Li-Yau-Hamilton Harnack inequality on compact Riemannian manifolds equipped with variants of the $(K, m)$-super Ricci flow. 
 In the literature, the Li-Yau Harnack inequality for heat equation $\partial_t u=\Delta u$ on compact Ricci flow has been studied by many authors. See \cite{Cao, CCLLN, Sun} and references therein.

\subsection{The commutator $[\partial_t,  L]f$}

Let $M$ be a compact manifold with a family of time dependent metrics $(g(t), t\in [0, T])$ and potentials $\phi(t)\in C^2(M)$, $t\in [0, T]$. Let $\partial_t g=2h$.

\begin{lemma} \label{LLLL} For any $f\in C^\infty(M)$, it holds
\begin{eqnarray*}
\partial_t |\nabla f|^2=-2h(\nabla f, \nabla f)+2\langle \nabla f, \nabla f_t\rangle,
\end{eqnarray*}
and
\begin{eqnarray*}
[\partial_t, L] f=-2\langle h, \nabla^2 f\rangle+2h(\nabla \phi, \nabla f)-\langle 2 {\rm div}h-\nabla {\rm Tr}_g h+\nabla \partial_t \phi, \nabla f\rangle.
\end{eqnarray*}
\end{lemma}
{\it Proof}. By direct calculation, cf. \cite{CLN, Sun},  we have

\begin{eqnarray*}
\partial_t \Delta_{g(t)} f=\Delta_{g(t)}\partial_t f-2\langle h, \nabla^2 f\rangle-2\langle {\rm div}h- {1\over 2}\nabla{\rm Tr}_g h, \nabla f\rangle,
\end{eqnarray*}
and
\begin{eqnarray*}
\partial_t \langle \nabla \phi, \nabla f\rangle=-\partial_t g(\nabla \phi, \nabla f)+\langle \nabla\phi_t, \nabla f\rangle+\langle \nabla \phi, \nabla f_t\rangle.
\end{eqnarray*}
Therefore
\begin{eqnarray*}
\partial_t L f&=&\partial_t \Delta_{g(t)} f-\partial_t \langle \nabla \phi, \nabla f\rangle\\
&=&\Delta_{g(t)} \partial_t f-2\langle h, \nabla^2 f\rangle-2\langle {\rm div}h- {1\over 2}\nabla{\rm Tr}_g h, \nabla f\rangle\\
& &\hskip1cm +2h(\nabla \phi, \nabla f)-\langle \nabla\phi_t, \nabla f\rangle-\langle \nabla \phi, \nabla f_t\rangle\\
&=&L\partial_t f-2\langle h, \nabla^2 f\rangle+2h(\nabla \phi, \nabla f)-\langle 2{\rm div}h-\nabla {\rm Tr}_g h+\nabla \phi_t, \nabla f\rangle.
\end{eqnarray*}
This finishes the proof. \hfill $\square$

\subsection{The Li-Yau Harnack inequality}

Let $u$ be a positive solution to the heat equation $\partial_t u=Lu$. Let $f=\log u$. Then
\begin{eqnarray*}
(L-\partial_t) f=-|\nabla f|^2.
\end{eqnarray*}
Let
\begin{eqnarray*}
F=t(|\nabla f|^2-\alpha f_t).
\end{eqnarray*}
We have

\begin{lemma}
\begin{eqnarray*}
(L-\partial_t) F=2t\left( |\nabla^2 f|^2+(Ric(L)+(1-\alpha)h)(\nabla f, \nabla f)\right)-2\langle \nabla f, \nabla F\rangle-t^{-1}F+\alpha t[\partial_t, L]f.
\end{eqnarray*}
\end{lemma}
{\it Proof}. By the Bochner formula and using
\begin{eqnarray*}
\partial_t |\nabla f|_{g(t)}^2=-\partial_t g(t)(\nabla f, \nabla f)+2\langle \nabla f, \nabla f_t\rangle_{g(t)},
\end{eqnarray*}
we have
\begin{eqnarray*}
LF&=&tL|\nabla f|^2-\alpha t Lf_t\\
&=&2t\left( |\nabla^2 f|^2+Ric(L)(\nabla f, \nabla f)+\langle \nabla f, \nabla L f\rangle \right) -\alpha t L\partial_t f\\
&=&2t\left( |\nabla^2 f|^2+Ric(L)(\nabla f, \nabla f)+\langle \nabla f, \nabla (f_t-|\nabla f|^2)\rangle \right) -\alpha t L\partial_t f\\
&=&2t\left( |\nabla^2 f|^2+Ric(L)(\nabla f, \nabla f)\right)-2\langle \nabla f, \nabla F\rangle+2(1-\alpha) t\langle \nabla f, \nabla f_t\rangle-\alpha t L\partial_t f\\
&=&2t\left( |\nabla^2 f|^2+Ric(L)(\nabla f, \nabla f)\right)-2\langle \nabla f, \nabla F\rangle\\
& &\hskip1cm +2t(1-\alpha) h(\nabla f, \nabla f)+(1-\alpha)t \partial_t |\nabla f|^2-\alpha t L\partial_t f.
\end{eqnarray*}
On the other hand
\begin{eqnarray*}
\partial_t F&=&(|\nabla f|^2-\alpha f_t)+t\partial_t |\nabla f|^2-\alpha t f_{tt}\\
&=&(|\nabla f|^2-\alpha f_t)+t\partial_t |\nabla f|^2-\alpha t \partial_t (Lf+|\nabla f|^2)\\
&=&(|\nabla f|^2-\alpha f_t)+(1-\alpha)t\partial_t |\nabla f|^2-\alpha t \partial_t Lf.
\end{eqnarray*}
Thus
\begin{eqnarray*}
(L-\partial_t) F=2t\left( |\nabla^2 f|^2+(Ric(L)+(1-\alpha)h)(\nabla f, \nabla f)\right)-2\langle \nabla f, \nabla F\rangle-t^{-1}F+\alpha t[\partial_t, L]f.
\end{eqnarray*}
By Lemma \ref{LLLL}, it holds
\begin{eqnarray*}
[\partial_t, L] f=-2 \left\langle h, \nabla^2 f\right\rangle+2h(\nabla \phi, \nabla f\rangle-\langle 2{\rm div} h-\nabla {\rm Tr}_g h+\nabla \phi_t, \nabla f\rangle.
\end{eqnarray*}
Thus
\begin{eqnarray*}
(L-\partial_t) F&=&2t\left|\nabla^2 f-{\alpha h\over 2}\right|^2-{t\alpha^2|h|^2\over 2} +2t(Ric(L)+(1-\alpha)h)(\nabla f, \nabla f)\\
& &-2\langle \nabla f, \nabla F\rangle-t^{-1}F+\alpha t S_1(\nabla f),
\end{eqnarray*}
where
\begin{eqnarray*}
S_1(\nabla f)=2h(\nabla \phi, \nabla f\rangle-\langle 2{\rm div} h-\nabla {\rm Tr}_g h+\nabla \phi_t, \nabla f\rangle.
\end{eqnarray*}
Using the Cauchy-Schwartz inequality, for all $\varepsilon>0$, we have $(a+b)^2\geq {a^2\over 1+\varepsilon}-{b^2\over \varepsilon}$. Thus
\begin{eqnarray*}
\left|\nabla^2 f-{\alpha h\over 2}\right|^2&\geq &{1\over n}\left|\Delta f-{\alpha {\rm Tr}h\over 2} \right|^2\\
&\geq& {|Lf|^2\over n(1+\varepsilon)}-{\left|\nabla \phi \cdot\nabla f-{\alpha {\rm Tr}h \over 2}\right|^2\over n\varepsilon}.
\end{eqnarray*}
Let $m:=n(1+\varepsilon)$. Then
\begin{eqnarray}
(L-\partial_t) F&\geq &{2t\over m}|Lf|^2-{2t\over m-n}\left|\nabla \phi \cdot\nabla f-{\alpha {\rm Tr}h \over 2}\right|^2 -{t\alpha^2|h|^2\over 2} +2t(Ric(L)+(1-\alpha)h)(\nabla f, \nabla f)\nonumber\\
& &-2\langle \nabla f, \nabla F\rangle-t^{-1}F+\alpha t S_1(\nabla f)\nonumber\\
&=&{2t\over m}|Lf|^2-{t\alpha^2 ({\rm Tr}h)^2\over 2(m-n)}-{t\alpha^2|h|^2\over 2} +2t(Ric_{m, n}(L)+(1-\alpha)h)(\nabla f, \nabla f)\nonumber\\
& &-2\langle \nabla f, \nabla F\rangle-t^{-1}F+\alpha t S_1(\nabla f)+{2\alpha t {\rm Tr} h\over m-n}\langle \nabla\phi, \nabla f\rangle.\label{FF1}
\end{eqnarray}
Let
\begin{eqnarray*}
S(\cdot)=S_1(\cdot)+{2 {\rm Tr} h\over m-n}\langle \nabla\phi, \cdot\rangle.
\end{eqnarray*}
Substituting $Lf=|\nabla f|^2-f_t={F\over \alpha t}+{\alpha-1\over \alpha}|\nabla f|^2$ into $(\ref{FF1})$, we have
\begin{eqnarray*}
(L-\partial_t) F
&\geq &{2t\over m}\left[{F\over \alpha t}+{\alpha-1\over \alpha}|\nabla f|^2\right]^2-{t\alpha^2\over 2}\left[{({\rm Tr}h)^2\over m-n}+|h|^2 \right]\\
& &+2t(Ric_{m, n}(L)+(1-\alpha)h)(\nabla f, \nabla f)-2\langle \nabla f, \nabla F\rangle-t^{-1}F+\alpha t S(\nabla f).
\end{eqnarray*}
Hence, as $F\geq 0$, and $\alpha\geq 1$, we have
\begin{eqnarray}
(L-\partial_t) F&\geq &{2F^2\over \alpha^2 m t}+{4(\alpha-1)\over
m\alpha^2}|\nabla f|^2F
+{2t(\alpha-1)^2\over m\alpha^2}|\nabla f|^4-{t\alpha^2\over 2}\left[{({\rm Tr}h)^2\over m-n}+|h|^2 \right]\nonumber\\
& &+2t(Ric_{m, n}(L)+(1-\alpha)h)(\nabla f, \nabla f)-2\langle
\nabla f, \nabla F\rangle-t^{-1}F+\alpha t S(\nabla
f).\nonumber \\
& &\hfill \label{ooooo}
\end{eqnarray}
Under the assumption
\begin{eqnarray}
Ric_{m, n}(L)+(1-\alpha)h\geq -K, \label{hh1}
\end{eqnarray}
and setting
\begin{eqnarray*}
A^2=\max\limits \left[ |h|^2+{({\rm Tr}h)^2\over m-n}\right],\ \ \ \ \ \ B=\max\limits |S|,
\end{eqnarray*}
we have
\begin{eqnarray*}
(L-\partial_t) F&\geq &{2F^2\over \alpha^2 m t}+{2t(\alpha-1)^2\over m\alpha^2}|\nabla f|^4-{t\alpha^2A^2\over 2}\\
& &-2Kt|\nabla f|^2-2\langle \nabla f, \nabla F\rangle-t^{-1}F-\alpha B t |\nabla f|.
\end{eqnarray*}
Using the inequality
\begin{eqnarray}
ax^4+bx^2+cx\geq -{(b-\gamma)^2\over 4a}-{c^2\over 4\gamma}, \label{xx1}
\end{eqnarray}
where $\gamma>0$ is any positive constant, we can derive that
\begin{eqnarray}
{2t(\alpha-1)^2\over m\alpha^2}|\nabla f|^4-2Kt|\nabla f|^2-\alpha B
t |\nabla f|\geq -{m\alpha^2 t(2K+\gamma)^2\over
8(\alpha-1)^2}-{\alpha^2 B^2 t\over 4\gamma}.\label{0}
\end{eqnarray}
Hence
\begin{eqnarray*}
(L - \partial_t)F
&\geq &  {2F^2\over \alpha^2 mt}-\frac{F}{t} - 2\langle \nabla f, \nabla F\rangle -{t\alpha^2A^2\over 2}\\
& &\ \ -{m\alpha^2 t(2K+\gamma)^2\over 8(\alpha-1)^2}-{\alpha^2 B^2 t\over 4\gamma}.
\end{eqnarray*}

Let $(x_0, t_0)$ be the point where $F$ achieves the maximum on $M\times [0, T]$.
Then $\nabla F(x_0, t_0)=0$, $\Delta F(x_0, t_0)\leq 0$ and $\partial_t F(x_0, t_0)\geq 0$. Therefore, at $(x_0, t_0)$,
\begin{eqnarray*}
(L-\partial_t)F\leq 0,
\end{eqnarray*}
i.e.,
\begin{eqnarray}
0\geq {2F^2\over \alpha^2 m t_0}-{F\over t_0} -{t_0\alpha^2A^2\over 2}-{m\alpha^2 t_0(2K+\gamma)^2\over 8(\alpha-1)^2}-{\alpha^2 B^2 t_0\over 4\gamma}.
\label{xx2}
\end{eqnarray}
Thus
\begin{eqnarray}
F\leq {m\alpha^2\over 4}\left[1+\sqrt{1+{t_0^2\over m}\left(4A^2+{(2K+\gamma)^2\over (\alpha-1)^2}+{2B^2\over \gamma}\right)}\right]. \label{xx3}
\end{eqnarray}
Note that when $B=0$, we can take $\gamma=0$ in $(\ref{xx1})$, $(\ref{xx2})$ and $(\ref{xx3})$, i.e.,
\begin{eqnarray}
F\leq  {m\alpha^2\over 4}\left[1+\sqrt{1+{t_0^2\over m}\left(4A^2+{4K^2\over (\alpha-1)^2}\right)}\right]. \label{xx4}
\end{eqnarray}
and if $K=0$, i.e., if
\begin{eqnarray}
Ric_{m, n}(L)+(1-\alpha)h\geq 0, \label{hh2}
\end{eqnarray}
we have
\begin{eqnarray}
F\leq {m\alpha^2 \over 4}\left[1+\sqrt{1+{4T^2A^2\over m}} \right]. \label{xx5}
\end{eqnarray}
In particular, when $A=B=0$, and $Ric_{m, n}(L)\geq 0$,  we can take $\alpha\rightarrow 1$ and recapture the Li-Yau Harnack inequality \cite{LY}
\begin{eqnarray}
|\nabla f|^2-f_t\leq {m\over 2t}. \label{xx6}
\end{eqnarray}

Therefore, we have proved the following Li-Yau Harnack inequality for positive solutions to the heat equation $\partial_t u=Lu$ on compact Riemannian manifolds equipped with  the backward $(K, m)$-super Ricci flows.

\begin{theorem}\label{LYHSRF-B}  Let $(M, g(t), t\in [0, T])$ be a compact Riemannian manifold with a family of time dependent metrics $g(t)$ and potentials $\phi(t)\in C^2(M)$, $t\in [0, T]$. Let $u$ be a positive solution to the heat equation $\partial_t u=Lu$. Let $\partial_t g=2h$ and $\alpha>1$. Suppose that
\begin{eqnarray}
{1\over 2}(1-\alpha)\partial_t g+Ric_{m, n}(L)\geq -Kg,\label{mmm1}
\end{eqnarray}
and assuming that $A^2=\max\limits \left[ |h|^2+{({\rm Tr}h)^2\over m-n}\right]<\infty$ and $B=\max\limits |S|<\infty$, where
\begin{eqnarray*}
S(\cdot)=2h(\nabla \phi, \cdot\rangle-\langle 2{\rm div} h-\nabla {\rm Tr}_g h+\nabla \phi_t, \cdot\rangle+{2 {\rm Tr} h\over m-n}\langle \nabla\phi, \cdot\rangle.
\end{eqnarray*}
Then for any $\gamma>0$ and for all $t\in (0, T]$, we have

\begin{eqnarray*}
{|\nabla u|^2\over u^2}-\alpha {\partial_t u\over u}\leq {m\alpha^2\over 4t}\left[1+\sqrt{1+{T^2\over m}\left(4A^2+{(2K+\gamma)^2\over (\alpha-1)^2}+{2B^2\over \gamma}\right)}\right].
\end{eqnarray*}
In the case $B=0$,  for all $t\in (0, T]$, we have
\begin{eqnarray*}
{|\nabla u|^2\over u^2}-\alpha {\partial_t u\over u}\leq {m\alpha^2 \over 4t}\left[1+\sqrt{1+{T^2\over m}\left(4A^2+{4K^2\over (\alpha-1)^2}\right)}\right].
\end{eqnarray*}
In particular, in the case $A=B=0$ and $Ric_{m, n}(L)\geq 0$, we have the Li-Yau Harnack inequality
\begin{eqnarray*}
{|\nabla u|^2\over u^2}-{\partial_t u\over u}\leq {m\over 2t}.
\end{eqnarray*}
\end{theorem}

\subsection{The Li-Yau-Hamilton Harnack inequality}

Let $u$ be a positive solution to the heat equation $\partial_t u=Lu$. Let $f=\log u$. Then $(\partial_t - L)f = |\nabla f|^{2}$. Let
\begin{eqnarray*}
F=te^{-2Kt}(e^{-2Kt}|\nabla f|^2-f_t)=te^{-4Kt}|\nabla f|^2-te^{-2Kt}f_t.
\end{eqnarray*}

In this section we prove the following Li-Yau-Hamilton Harnack inequality on a variant of the $(K, m)$-super Ricci flow on compact manifolds.

\begin{theorem}\label{LYHHSRF}  Let $(M, g(t), t\in [0, T])$ be a compact Riemannian manifold with a family of time dependent metrics $g(t)$ and potentials $\phi(t)\in C^2(M)$, $t\in [0, T]$. Let $u$ be a positive solution to the heat equation $\partial_t u=Lu$. Suppose that $\partial_t g=2h$  satisfies 
\begin{eqnarray}
e^{-4Kt} (h+Ric_{m, n}(L)+Kg)-e^{-2Kt}h\geq  \alpha_K(t)g,
\end{eqnarray}
and
\begin{eqnarray*}
A^2=\max\limits \left[ |h|^2+{({\rm Tr}h)^2\over m-n}\right]<\infty,\ \ \ \ \ \ B=\max\limits |S|<\infty,
\end{eqnarray*}
where
\begin{eqnarray*}
S(\cdot)=\left\langle {2{\rm Tr}h\over m-n} \nabla \phi-2 {\rm div}h+\nabla {\rm Tr}_g h-\nabla \partial_t \phi, \cdot\right\rangle +2h(\nabla\phi, \cdot).
\end{eqnarray*}
Then, for any $\gamma>0$ and $t\in [0, T]$, we have
\begin{eqnarray*}
{|\nabla u|^2\over u^2}-e^{2Kt}{\partial_t u\over u}  \leq  {me^{4Kt}\over 2t} \left[1+\sqrt{{A^2T^2\over m}+\max\limits_{t\in [0, T]}{t^2(2\alpha_K(t)-\gamma)^2\over 4e^{-4Kt}(1-e^{-2Kt})^2}+\max\limits_{t\in [0, T]}{t^2e^{-4Kt}B^2\over 2m\gamma} } \right].
\end{eqnarray*}
In the case $B=0$, we have
\begin{eqnarray*}
{|\nabla u|^2\over u^2}-e^{2Kt}{\partial_t u\over u}  \leq  {me^{4Kt}\over 2t} \left[1+\sqrt{{A^2T^2\over m}+\max\limits_{t\in [0, T]}{t^2\alpha^2_K(t)\over e^{-4Kt}(1-e^{-2Kt})^2} } \right].
\end{eqnarray*}
and if $\alpha_K(t)=0$, i.e., if
\begin{eqnarray}
e^{-4Kt} (h+Ric_{m, n}(L)+K)-e^{-2Kt}h\geq  0,
\end{eqnarray}
we have
\begin{eqnarray*}
{|\nabla u|^2\over u^2}-e^{2Kt}{\partial_t u\over u}  \leq  {me^{4Kt}\over 2t} \left[1+{TA\over \sqrt{m}} \right].
\end{eqnarray*}
In particular, when $A=B=0$, and $Ric_{m, n}(L)\geq -K$,  we recapture Hamilton's Harnack inequality \cite{H1}
\begin{eqnarray*}
{|\nabla u|^2\over u^2}-e^{2Kt}{\partial_t u\over u}  \leq {me^{4Kt}\over 2t}.
\end{eqnarray*}
\end{theorem}

{\it Proof}.  If $F\leq 0$ on $[0, T]\times M$, we have
${|\nabla u|^2\over u^2}-e^{2Kt} {\partial_t u\over u}\leq 0$.
In this case, the Li-Yau-Hamilton Harnack inequality automatically holds. Thus, in order to prove Theorem \ref{LYHHSRF}, we need only to consider the case $F\geq 0$. Notice that
\begin{eqnarray*}
L|\nabla f|^{2} &=& 2|Hess f|^{2} + 2\langle \nabla f, \nabla L f\rangle + 2Ric(L)(\nabla f, \nabla f),\\
\partial_t |\nabla f|^{2} &=& 2\langle \nabla \partial_t f, \nabla f\rangle - 2h(\nabla f, \nabla f),\\
\langle \nabla F, \nabla f\rangle &=& 2te^{-4Kt}(1- e^{2Kt})\nabla^{2}f(\nabla f, \nabla f) - te^{-2Kt}\langle \nabla L f, \nabla f\rangle\\
&=& te^{-2Kt}[2(e^{-2Kt}-1)\nabla^{2}f(\nabla f, \nabla f) - \langle \nabla L f, \nabla f\rangle].
\end{eqnarray*}
Then
\begin{eqnarray*}
(\partial_t - L)F &=& \partial_t(te^{-2Kt}(e^{-2Kt}|\nabla f|^2-f_t)) - te^{-2Kt}(e^{-2Kt}L|\nabla f|^2 - Lf_t)\\
&=& e^{-2Kt}(1 - 2Kt)(e^{-2Kt}|\nabla f|^2-f_t) + te^{-2Kt}[e^{-2Kt}\partial_t|\nabla f|^2-2Ke^{-2Kt}|\nabla f|^{2}-f_{tt})\\
& &- te^{-2Kt}(e^{-2Kt}L|\nabla f|^2 - Lf_t)\\
&=& \frac{1-2Kt}{t}F + te^{-2Kt}[e^{-2Kt}(\partial_t -L)|\nabla f|^{2} - 2Ke^{-2Kt}|\nabla f|^{2} - (\partial_t - L)f_t]\\
&=& \frac{1-2Kt}{t}F + te^{-2Kt}[e^{-2Kt}(\partial_t -L)|\nabla f|^{2} - 2Ke^{-2Kt}|\nabla f|^{2} - \partial_t |\nabla f|^{2} - [\partial_t, L]f]\\
&=& \frac{1-2Kt}{t}F + te^{-2Kt}[(e^{-2Kt}-1)\partial_t|\nabla f|^{2} -e^{-2Kt}L|\nabla f|^{2} - 2Ke^{-2Kt}|\nabla f|^{2} - [\partial_t, L]f]\\
&=& \frac{1-2Kt}{t}F + te^{-2Kt}[(e^{-2Kt}-1)(2\langle \nabla \partial_t f, \nabla f\rangle - 2h(\nabla f, \nabla f)) \\
& &-e^{-2Kt}(2|Hess f|^{2} + 2\langle \nabla f, \nabla L f\rangle + 2Ric(L)(\nabla f, \nabla f)) - 2Ke^{-2Kt}|\nabla f|^{2} - [\partial_t, L]f]\\
&=& \frac{1-2Kt}{t}F + te^{-2Kt}[(e^{-2Kt}-1)(2\langle \nabla L f, \nabla f \rangle + 2\langle \nabla |\nabla f|^{2}, \nabla f \rangle - 2h(\nabla f, \nabla f)) \\
& & -e^{-2Kt}(2|Hess f|^{2} + 2\langle \nabla f, \nabla L f\rangle + 2Ric(L)(\nabla f, \nabla f)) - 2Ke^{-2Kt}|\nabla f|^{2} - [\partial_t, L]f]\\
&=& \frac{1-2Kt}{t}F + te^{-2Kt}[(e^{-2Kt}-1)(4\nabla^{2}f(\nabla f, \nabla f) - 2h(\nabla f, \nabla f)) - 2\langle \nabla f, \nabla L f\rangle\\
& & -e^{-2Kt}(2|Hess f|^{2}+ 2(Ric(L)+K)(\nabla f, \nabla f)) - [\partial_t, L]f]\\
&=& \frac{1-2Kt}{t}F + 2\langle \nabla F, \nabla f\rangle + te^{-2Kt}[-e^{-2Kt}(2|Hess f|^{2}+ 2(h+Ric(L)+K)(\nabla f, \nabla f)) \\
& & + 2h(\nabla f, \nabla f) - [\partial_t, L]f].
\end{eqnarray*}
This yields
\begin{eqnarray*}
(L - \partial_t)F &=& \frac{2Kt-1}{t}F - 2\langle \nabla F, \nabla f\rangle + te^{-2Kt}[e^{-2Kt}(2|Hess f|^{2}+ 2(h+Ric(L)+K)(\nabla f, \nabla f)) \\
& & -2h(\nabla f, \nabla f) - 2\langle h, \nabla^2 f\rangle + 2h(\nabla \phi, \nabla f) - \langle 2 {\rm div}h-\nabla {\rm Tr}_g h+\nabla \partial_t \phi, \nabla f\rangle]\\
&=& \frac{2Kt-1}{t}F - 2\langle \nabla F, \nabla f\rangle + 2te^{-4Kt}\left|\nabla^2 f - \frac{e^{2Kt}}{2}h\right|^{2} - \frac{t}{2}|h|^{2} \\
& &\ +2te^{-4Kt} (h+Ric(L)+K)(\nabla f, \nabla f))-2te^{-2Kt}h(\nabla f, \nabla f) \\
& &\  +te^{-2Kt}(2h(\nabla \phi, \nabla f) - \langle 2 {\rm div}h-\nabla {\rm Tr}_g h+\nabla \partial_t \phi, \nabla f\rangle).
\end{eqnarray*}  
Note that
\begin{eqnarray*}
\left|\nabla^2 f - \frac{e^{2Kt}}{2}h\right|^{2} \geq \frac{1}{n}\left|\Delta f-{e^{2Kt}\over 2}{\rm Tr}h\right|^{2} \\
\end{eqnarray*}
Using the elementary inequality   $(a+b)^2\geq {a^2\over 1+\varepsilon}-{b^2\over \varepsilon}$ with $a=Lf$, $b=\nabla \phi\cdot\nabla f-{e^{2Kt}\over 2}{\rm Tr}h$,  we have
\begin{eqnarray*}
\left|\Delta f-{e^{2Kt}\over 2}{\rm Tr}h\right|^{2}
 \geq \frac{\ |Lf|^{2} }{1+\varepsilon}-{1\over \varepsilon}\left|\nabla \phi\cdot\nabla f-{e^{2Kt}\over 2}{\rm Tr}h\right|^2
\end{eqnarray*}
Taking $\varepsilon ={m-n\over n}$, we obtain
\begin{eqnarray*}
\left|\nabla^2 f - \frac{e^{2Kt}}{2}h\right|^{2} \geq \frac{|Lf|^{2}}{m} -{1\over m-n}\left|\nabla \phi\cdot\nabla f-{e^{2Kt}\over 2}{\rm Tr}h\right|^2.
\end{eqnarray*}
This yields
\begin{eqnarray*}
(L - \partial_t)F &\geq &
\frac{2Kt-1}{t}F - 2\langle \nabla F, \nabla f\rangle + {2te^{-4Kt} \over m}|Lf|^{2}-{2te^{-4Kt}\over m-n}\left|\nabla \phi\cdot\nabla f-{e^{2Kt}\over 2}{\rm Tr}h\right|^2 \\
& &\ \  +2te^{-4Kt} (h+Ric(L)+K)(\nabla f, \nabla f))-2te^{-2Kt}h(\nabla f, \nabla f) - \frac{t}{2}|h|^{2}\\
& &\ \  +te^{-2Kt}(2h(\nabla \phi, \nabla f) - \langle 2 {\rm div}h-\nabla {\rm Tr}_g h+\nabla \partial_t \phi, \nabla f\rangle)\\
&=&\frac{2Kt-1}{t}F - 2\langle \nabla F, \nabla f\rangle + {2te^{-4Kt} \over m}|Lf|^{2}-{t\over 2}\left[{({\rm Tr}h)^2\over m-n}+|h|^2\right]\\
& &\ +2te^{-4Kt} (h+Ric_{m, n}(L)+K)(\nabla f, \nabla f))-2te^{-2Kt}h(\nabla f, \nabla f) \\
& &\  +te^{-2Kt}\left({2{\rm Tr}h\over m-n}\langle \nabla \phi, \nabla f\rangle +2h(\nabla \phi, \nabla f) - \langle 2 {\rm div}h-\nabla {\rm Tr}_g h+\nabla \partial_t \phi, \nabla f\rangle\right).
\end{eqnarray*}
Substituting $Lf=(e^{-2Kt}-1)|\nabla f|^{2} - \frac{e^{2Kt}}{t}F$ into the above inequality, we get
\begin{eqnarray*}
(L - \partial_t)F
&\geq & \frac{2Kt-1}{t}F - 2\langle \nabla F, \nabla f\rangle -{t\over 2}\left[{({\rm Tr}h)^2\over m-n}+|h|^2\right]\\
& &\ + {2te^{-4Kt} \over m}\left| (e^{-2Kt}-1)|\nabla f|^{2} - \frac{e^{2Kt}}{t}F \right|^{2}\\
& &\  +2te^{-4Kt} (h+Ric_{m, n}(L)+K)(\nabla f, \nabla f))-2te^{-2Kt}h(\nabla f, \nabla f) \\
& &\  +te^{-2Kt}\left({2{\rm Tr}h\over m-n}\langle \nabla \phi, \nabla f\rangle +2h(\nabla \phi, \nabla f) - \langle 2 {\rm div}h-\nabla {\rm Tr}_g h+\nabla \partial_t \phi, \nabla f\rangle\right).
\end{eqnarray*}
Under the assumption
\begin{eqnarray}
e^{-4Kt} (h+Ric_{m, n}(L)+K)-e^{-2Kt}h\geq  \alpha_K(t),
\end{eqnarray}
where $\alpha_K(t)$ is a function of $t$ and $K$, we have (using the assumption $F\geq 0$ and $K\geq 0$)
\begin{eqnarray}
(L - \partial_t)F
&\geq & \frac{2Kt-1}{t}F - 2\langle \nabla F, \nabla f\rangle -{t\over 2}\left[{({\rm Tr}h)^2\over m-n}+|h|^2\right]  + {2F^2\over mt}\nonumber\\
& &+{4e^{-4Kt} (1-e^{-2Kt})^2|\nabla f|^{2}F\over m}+ {2te^{-4Kt} (1-e^{-2Kt})^2|\nabla f|^{4}\over m} +2t\alpha_K(t)|\nabla f|^2\nonumber\\
& &+te^{-2Kt}\left({2{\rm Tr}h\over m-n}\langle \nabla \phi, \nabla f\rangle +2h(\nabla \phi, \nabla f) - \langle 2 {\rm div}h-\nabla {\rm Tr}_g h+\nabla \partial_t \phi, \nabla f\rangle\right).\nonumber\\
& &\hskip4cm \label{GGGG1}
\end{eqnarray}
Set
\begin{eqnarray}
S(\cdot)=\left\langle {2{\rm Tr}h\over m-n} \nabla \phi-2 {\rm div}h-\nabla {\rm Tr}_g h+\nabla \partial_t \phi, \cdot\right\rangle +2h(\nabla\phi, \cdot), \label{SSS}
\end{eqnarray}
and
\begin{eqnarray*}
A^2=\max\limits \left[ |h|^2+{({\rm Tr}h)^2\over m-n}\right],\ \ \ \ \ \ B=\max\limits |S|.
\end{eqnarray*}
Then
\begin{eqnarray*}
(L - \partial_t)F
&\geq & \frac{2Kt-1}{t}F - 2\langle \nabla F, \nabla f\rangle -{tA^2\over 2} + {2F^2\over mt}\\
& &\ + {2te^{-4Kt} (1-e^{-2Kt})^2|\nabla f|^{4}\over m} +2t\alpha_K(t)|\nabla f|^2-te^{-2Kt}B|\nabla f|.
\end{eqnarray*}
Using the inequality
\begin{eqnarray}
ax^4+bx^2+cx\geq -{(b-\gamma)^2\over 4a}-{c^2\over 4\gamma}, \label{xxx1}
\end{eqnarray}
where $\gamma>0$ is any positive constant, we can derive that
\begin{eqnarray*}
& &{2te^{-4Kt} (1-e^{-2Kt})^2|\nabla f|^{4}\over m} +2t\alpha_K(t)|\nabla f|^2-te^{-2Kt}B|\nabla f|\\
& &\hskip2cm  \geq -{mt(2\alpha_K(t)-\gamma)^2\over 8e^{-4Kt}(1-e^{-2Kt})^2}-{te^{-4Kt}B^2\over 4\gamma}.
\end{eqnarray*}
Thus
\begin{eqnarray}
(L - \partial_t)F
&\geq & \frac{2Kt-1}{t}F - 2\langle \nabla F, \nabla f\rangle -{tA^2\over 2} + {2F^2\over mt}\nonumber\\
& &\ \ -{mt(2\alpha_K(t)-\gamma)^2\over 8e^{-4Kt}(1-e^{-2Kt})^2}-{te^{-4Kt}B^2\over 4\gamma}.\label{GGGG2}
\end{eqnarray}
Suppose at $(x_0, t_0) \in M \times [0, T]$, $F$ achieves its maximum. Then $F(x_0, t_0) \geq 0$, and at $(x_0, t_0)$ we have
$$
\partial_t F \geq 0, \ \ \ \Delta F \leq 0, \ \ \ \nabla F =0.
$$
Thus
$$
(L - \partial_t)F \leq 0.
$$
Multiplying $t_0$ on the both sides of the last inequality, we have
\begin{eqnarray}
0\geq  (2Kt_0-1) F+{2F^2\over m}-{t_0^2A^2\over 2}-{mt_0^2(2\alpha_K(t_0)-\gamma)^2\over 8e^{-4Kt_0}(1-e^{-2Kt_0})^2}-{t_0^2e^{-4Kt_0}B^2\over 4\gamma},
\end{eqnarray}
Thus we obtain the following Li-Yau-Harmilton Harnack inequality on super Ricci flow
\begin{eqnarray}
F&\leq & {m \over 4}\left[(1-2Kt_0)+\sqrt{(1-2Kt_0)^2+ {8t_0^2\over m}\left({A^2\over 2}+{m(2\alpha_K(t_0)-\gamma)^2\over 8e^{-4Kt_0}(1-e^{-2Kt_0})^2}+{e^{-4Kt_0}B^2\over 4\gamma}\right)  } \right]\nonumber\\
&\leq& {m\over 2}\left[(1-2Kt_0)^{+}+\sqrt{{2\over m}\left({T^2A^2\over 2}+{mt_0^2(2\alpha_K(t_0)-\gamma)^2\over 8e^{-4Kt_0}(1-e^{-2Kt_0})^2}+{t_0^2e^{-4Kt_0}B^2\over 4\gamma}\right)  } \right]\nonumber\\
&\leq& {m\over 2}\left[1+\sqrt{{A^2T^2\over m}+\max\limits_{t\in [0, T]}{t^2(2\alpha_K(t)-\gamma)^2\over 4e^{-4Kt}(1-e^{-2Kt})^2}+\max\limits_{t\in [0, T]}{t^2e^{-4Kt}B^2\over 2m\gamma} } \right].\label{xxx2}
\end{eqnarray}
Note that when $B=0$, we can take $\gamma=0$ in $(\ref{xxx1})$ and in $(\ref{xxx2})$, i.e.,
\begin{eqnarray}
F\leq {m\over 2}\left[1+\sqrt{{A^2T^2\over m}+\max\limits_{t\in [0, T]}{ t^2\alpha_K^2(t)e^{4Kt} \over (1-e^{-2Kt})^2}} \right], \label{xxx3}
\end{eqnarray}
and if $\alpha_K(t)=0$, i.e., if
\begin{eqnarray}
e^{-4Kt} (h+Ric_{m, n}(L)+K)-e^{-2Kt}h\geq  0,
\end{eqnarray}
we have
\begin{eqnarray}
F\leq {m\over 2}\left[1+{TA\over \sqrt{m}} \right]. \label{xxx4}
\end{eqnarray}
In particular, when $A=B=0$, and $Ric_{m, n}(L)\geq -K$,  we recapture Hamilton's Harnack inequality \cite{H1}
\begin{eqnarray}
F\leq {m\over 2}. \label{xxx5}
\end{eqnarray}

\subsection{Hamilton's second order estimates for time dependent Witten Laplacian}

Let $(M, g(t), \phi(t), t\in [0, T])$  be a compact Riemannian manifold equipped with a family of time dependent metrics $g(t)$ and potentials $\phi(t)$ , $t\in [0, T]$. Let 
\begin{eqnarray*}
\partial_t g=2h.
\end{eqnarray*}
Let $u$ be a positive solution to the heat equation $\partial_t u=L u$ associated with the time dependent Witten Laplacian $L=\Delta_{g(t)}-\nabla \phi(t)\cdot \nabla$. 
Let $P:=(\partial_t-L-2\nabla \log u\cdot \nabla )$, and let $\psi: [0, T]\rightarrow [0, \infty)$ be a $C^1$-function. Set
\begin{eqnarray*}
F(x, t)=\psi(t) \left({L u\over u}+{|\nabla u|^2\over u^2}\right)-\left(m+4\log \left({A\over u}\right)\right).
\end{eqnarray*}
By the Bochner formula, we have
\begin{eqnarray*}
PL\log u&=&\partial_t L_{t}\log u+L(\partial_t-L)\log u-2\nabla \log u\cdot \nabla L\log u\\
&=&\partial_t L_{t}\log u+L|\nabla \log u|^2-2\nabla \log u\cdot \nabla L\log u\\
&=&\partial_t L_{t}\log u+2|\nabla^2\log u|^2+2Ric(L)(\nabla \log u, \nabla u).
\end{eqnarray*}
On the other hand
\begin{eqnarray*}
P |\nabla\log u|^2&=&-\partial_t g(\nabla \log u, \nabla \log u)+2\nabla \log u\cdot \nabla\partial_t \log u-L|\nabla\log u|^2-2\nabla  \log u \cdot \nabla |\nabla \log u|^2\\
&=&-\partial_t g(\nabla \log u, \nabla \log u)+2\nabla \log u\cdot \nabla L\log u-L|\nabla\log u|^2\\
&=&-\left(\partial_t g+2Ric(L)\right)(\nabla\log u, \nabla \log u)-2|\nabla^2\log u|^2.
\end{eqnarray*}
Combining the above two formulas together we have 
\begin{eqnarray*}
PF &=&\psi(t)P \left({Lu\over u}+{|\nabla u|^2\over u^2}\right)+\psi'(t)\left( {L u\over u}+{|\nabla u|^2\over u^2}\right)-4P \log \left({A\over u}\right)\\
&=&\psi(t)\left(\partial_t L_t \log u-2|\nabla^2\log u|^2-2(\partial_t g+Ric(L))(\nabla \log u, \nabla \log u)\right)\\
& &\ \ \ \ +\psi'(t)(L\log u+2|\nabla\log u|^2)-4|\nabla\log u|^2.
\end{eqnarray*} 
By Lemma \ref{LLLL}, we have 
\begin{eqnarray*}
\partial_t L_t \log u=-\langle \partial_t g, \nabla^2 \log u\rangle-S_1(\nabla \log u),
\end{eqnarray*}
where
\begin{eqnarray*}
S_1(\nabla \log u)=\langle {\rm div}\partial_t g-{1\over 2}\nabla {\rm Tr}\partial_t g-\nabla \partial_t \phi, \nabla \log u\rangle+\partial_t g(\nabla \phi, \nabla \log u).
\end{eqnarray*}
Thus
\begin{eqnarray*}
PF&= &-2\psi (|\nabla^2\log u|^2+\langle h, \nabla^2\log u\rangle) -2\psi (\partial_t g+Ric(L))(\nabla \log u, \nabla \log u)\\
& &\ \ \ \ +\psi' L \log u+2(\psi' -2)|\nabla \log u|^2-\psi S_1(\nabla \log u)\\
&= &-2\psi \left[|\nabla^2 \log u+{h\over 2}|^2-{|h|^2\over 4}\right]+\psi' L\log u+2(\psi'-(\partial_t g+Ric(L))\psi-2 )|\nabla \log u|^2-\psi S_1(\nabla \log u)\\
&\leq&-{2\psi \over n}|\Delta \log u+{{\rm Tr} h\over 2}|^2+{\psi |h|^2\over 2}+\psi' L \log u+2(\psi' -(\partial_t g+Ric(L))\psi-2)|\nabla \log u|^2-\psi S_1(\nabla \log u).
\end{eqnarray*}
Using the inequality $(a+b)^2\geq {a^2\over 1+\varepsilon}-{b^2\over \varepsilon}$, and taking $\varepsilon={m-n\over n}$ for any $m\geq n$, we have
\begin{eqnarray*}
{1\over n}|\Delta \log u+{{\rm Tr} h\over 2}|^2\geq {|L\log u|^2\over m}-{|\nabla\phi\cdot \nabla \log u +{{\rm Tr} h\over 2}|^2\over m-n}.
\end{eqnarray*}
Therefore
\begin{eqnarray*}
PF&\leq&-{2\psi \over m}|L \log u|^2+\psi' L \log u+{2\psi |\nabla\phi\cdot \nabla \log u+{{\rm Tr} h\over 2}|^2 \over m-n}+{\psi |h|^2\over 2}\\
& &\ \ \ \ +2(\psi'+Ric(L)\psi-2)|\nabla \log u|^2-\psi S_1(\nabla \log u)\\
&=& -{2\psi \over m}|L \log u|^2+\psi' L \log u+ {\psi \over 2}\left[{ |{\rm Tr} h|^2 \over m-n} +|h|^2\right]\\
& &\ \ \ \ +2(\psi'-(\partial_t g+Ric_{m, n}(L))\psi-2)|\nabla \log u|^2-\psi S(\nabla \log u),
\end{eqnarray*}
where 
\begin{eqnarray*}
S(\nabla \log u)=S_1(\nabla \log u)-{2{\rm Tr} h\over m-n}\nabla \phi\cdot\nabla \log u.
\end{eqnarray*}
Suppose that
\begin{eqnarray*}
\partial_t g+Ric_{m, n}(L)\geq -Kg.
\end{eqnarray*}
and denote $B=\max\limits\{ |S(v)|: v\in T _{\cdot}M, |v|=1\}$. We have 
\begin{eqnarray*}
PF&\leq& -{2\psi \over m}|L \log u|^2+\psi' L \log u+ {\psi \over 2}\left[{ |{\rm Tr} h|^2 \over m-n} +|h|^2\right]\\
& &\ \ \ \ +2(\psi'+K\psi-2)|\nabla \log u|^2+2\alpha |\nabla \log u|^2+{\psi^2B^2\over 8\alpha},
\end{eqnarray*}
where $\alpha$ is any constant with $\alpha\in (0, 1)$. Taking $\psi(t)=(1-\alpha){1-e^{-Kt}\over K}$, then 
\begin{eqnarray*}
PF&\leq&-{2\psi \over m}|L \log u|^2+\psi' L \log u-2|\nabla \log u|^2+{\psi^2 B^2\over 8\alpha}+{\psi \over 2}\left[{|{\rm Tr} h|^2 \over m-n}+|h|^2\right]\\
&\leq&-{2\psi \over m}|L \log u|^2+\psi' L \log u-2|\nabla \log u|^2+C,
\end{eqnarray*}
where 
\begin{eqnarray*}
C={(1-\alpha)^2B^2\over 8\alpha K^2}+{1-\alpha\over 2K}\left[{|{\rm Tr} h|^2 \over m-n}+|h|^2\right].
\end{eqnarray*}

Let $Q=-{2\psi \over m}|L \log u|^2+\psi' L \log u-2|\nabla \log u|^2+C$. Whenever $F\geq 0$, we have
\begin{eqnarray*}
\psi (L \log u+2|\nabla \log u|^2) \geq m+4\log \left({A\over u}\right)\geq m,
\end{eqnarray*}
which yields either $\psi L \log u\geq {m\over 2}$ or $2\psi|\nabla \log u|^2\geq {m\over 2}$. In the case $\psi L \log u\geq {m\over 2}$, we have
\begin{eqnarray*}
Q\leq (\psi'-1)L\log u-2|\nabla\log u|^2+C\leq C,
\end{eqnarray*} 
and in the case $2\psi|\nabla \log u|^2\geq {m\over 2}$,  we have
\begin{eqnarray*}
Q&\leq& -{2\psi\over m}|L \log u|^2+\psi' L \log u-{m\over 2\psi}+C\\
&\leq &{m\over 2\psi}\left({\psi^{'2}\over 4}-1\right)+C\\
&\leq& C.
\end{eqnarray*}
Thus, whenever $F\geq 0$,  we have
\begin{eqnarray*}
PF\leq Q\leq C.
\end{eqnarray*}
This yields that at any point where  $F\geq 0$ we have
\begin{eqnarray*}
P(F-Ct)\leq 0.
\end{eqnarray*}
Note that $F\leq 0$ at time $t=0$. By the maximum principle, we can derive that
\begin{eqnarray*}
F\leq Ct, \ \ \ \forall \ t\in [0, T].
\end{eqnarray*}

Therefore  we have proved the following Hamilton second order estimate for positive solutions to the heat equation $\partial_t u=L u$ on compact Riemannian manifolds equipped with (a variant of) the $(K, m)$-super Ricci flow.

\begin{theorem}\label{HHH} Let $M$ be a compact Riemannian manifold with a family of Riemannian metrics $(g(t)$ and potentials $\phi(t)$, $t\in [0, T]$. Suppose that for some constants $m>n$ and $K\in \mathbb{R}$,  
\begin{eqnarray*}
\partial_t g+Ric_{m, n}(L)\geq -K g, \ \ \ \forall \ t\in [0, T]. 
\end{eqnarray*}
 Let $u$ be a positive solution to the heat equation $\partial_t u=L u$. Let $A=\max\limits\{u(x, t): x\in M, t\in [0, T]\}$. Then for any $\alpha\in (0, 1)$, we have 
\begin{eqnarray*}
{L u\over u}+{|\nabla u|^2\over u^2}\leq {K\over (1-\alpha)(1-e^{-Kt})}\left[{m\over 2}+4\log\left({A\over u}\right)+Ct\right],
\end{eqnarray*}
where
\begin{eqnarray*}
C={(1-\alpha)^2B^2\over 8\alpha K^2}+{1-\alpha\over 2K}\max\limits\left[{|{\rm Tr} h|^2 \over m-n}+|h|^2\right],
\end{eqnarray*}
with $B=\max\limits\{ |S(v)|: v\in T _{x}M, |v|=1, x\in M, t\in [0, T]\}$, where 
\begin{eqnarray*}
S(v)=\left\langle {2{\rm Tr}h\over m-n} \nabla \phi-2 {\rm div}h+\nabla {\rm Tr}_g h-\nabla \partial_t \phi, v\right\rangle +2h(\nabla\phi, v).
\end{eqnarray*}
In the case  where $B=0$,  we can take $\alpha=0$, i.e., 
\begin{eqnarray*}
{L u\over u}+{|\nabla u|^2\over u^2}\leq {K\over 1-e^{-Kt}}\left[{m\over 2}+4\log\left({A\over u}\right)+Ct\right],
\end{eqnarray*}
where 
\begin{eqnarray*}
C={1\over 2K}\max\limits\left[{|{\rm Tr} h|^2 \over m-n}+|h|^2\right].
\end{eqnarray*}
In the case where $g$ and $\phi$ are time independent, and $Ric_{m, n}(L)\geq -Kg$, we have 
\begin{eqnarray*}
{L u\over u}+{|\nabla u|^2\over u^2}\leq {K\over 1-e^{-Kt}}\left[{m\over 2}+4\log\left({A\over u}\right)\right].
\end{eqnarray*}
\end{theorem}

In particular, when $L=\Delta$ and $\phi(t)\equiv 0$, we have the following 

\begin{theorem}\label{HHH-H} Let $M$ be a compact Riemannian manifold with a family of Riemannian metrics $(g(t)$, $t\in [0, T]$. Suppose that for some constant $K\in \mathbb{R}$,  
\begin{eqnarray*}
\partial_t g+Ric\geq -K g, \ \ \ \forall \ t\in [0, T]. 
\end{eqnarray*}
 Let $u$ be a positive solution to the heat equation $\partial_t u=\Delta u$. Let $A=\max\limits\{u(x, t): x\in M, t\in [0, T]\}$. Then for any $\alpha\in (0, 1)$ and any $m>n$, we have 
\begin{eqnarray*}
{\Delta u\over u}+{|\nabla u|^2\over u^2}\leq {K\over (1-\alpha)(1-e^{-Kt})}\left[{m\over 2}+4\log\left({A\over u}\right)+Ct\right],
\end{eqnarray*}
where
\begin{eqnarray*}
C={(1-\alpha)^2B^2\over 8\alpha K^2}+{1-\alpha\over 2K}\max\limits\left[{|{\rm Tr} h|^2 \over m-n}+|h|^2\right],
\end{eqnarray*}
with $B=\max\limits\{ |S(v)|: v\in T _{x}M, |v|=1, x\in M, t\in [0, T]\}$ for $S(v)=\langle {\rm div}\partial_t g-{1\over 2}\nabla {\rm Tr}\partial_t g, v\rangle$.  In the case where $g$ is the Ricci flow, i.e., $\partial_t g=-Ric$, we 
have $B=0$, and we can take $\alpha=0$, hence
\begin{eqnarray*}
{\Delta u\over u}+{|\nabla u|^2\over u^2}\leq {K\over 1-e^{-Kt}}\left[{m\over 2}+4\log\left({A\over u}\right)+Ct\right],
\end{eqnarray*}
where 
\begin{eqnarray*}
C={1\over 2K}\max\limits\left[{|{\rm Tr} h|^2 \over m-n}+|h|^2\right].
\end{eqnarray*}
In the case where $g$ and $\phi$ are time independent, and $Ric\geq -Kg$, we can take $m=n$ and we have Hamilton's second order estimate 
\begin{eqnarray*}
{\Delta u\over u}+{|\nabla u|^2\over u^2}\leq {K\over 1-e^{-Kt}}\left[{n\over 2}+4\log\left({A\over u}\right)\right].
\end{eqnarray*}
\end{theorem}

\section{The Li-Yau and the Li-Yau-Hamilton Harnack inequalities on complete super Ricci flow}

In this section we prove the Li-Yau Harnack inequality and the Li-Yau-Hamilton Harnack inequality on complete Riemannian manifolds equipped with variants of the $(K, m)$-super Ricci flow.  In the literature, the Li-Yau Harnack inequality for heat equation $\partial_t u=\Delta u$ on complete Ricci flow has been studied by many authors. See \cite{Cao, CCLLN, Sun} and references therein.

Similarly to \cite{Li05}, let $\eta$ be a $C^2$-function on $[0,―infty)$ such that $\eta=1$ on $[0, 1]$ and $\eta=0$ on $[2,―infty)$, with $-C_1\eta^{1/2}(r)\leq \eta'(r)\leq 0$, and $\eta^{''}(r)\geq C_2$, where $C_, C_2>0$ are two constants. Let $\rho(x)=d(o, x)$ and define $\psi(x)=\eta({\rho(x)/R})$.

Let $Q_{2R, T}=\{(x, t)\in M\times [0, T]: d(x, x_0, t)\leq 2R, t\in [0, T]\}$. Let $\eta\in C^2([0, \infty), [0, 1]) $ be such that $\eta(r)=1$ on $[0, 1]$, $\eta=0$ on $[2, \infty)$, $0\leq \eta \leq 1$ on $[1, 2]$, $\eta'(r)\leq 0$, $\eta''(r)\geq -C$ and $|\eta'(r)|^2\leq C\eta(r)$, where $C$ is a positive constant. Define
\begin{eqnarray*}
\psi(x, t)=\psi(d(x, x_0, t))=\eta\left({d(x, x_0, t)\over R}\right)=\eta\left({\rho(x, t)\over R}\right),
\end{eqnarray*}
where $\rho(x, t)=d(x, x_0, t)$ denotes the geodesic distance between $x$ and $x_0$ on $(M, g(t))$.

We need the following lemma.

\begin{lemma}\label{comparison2} Let $M$ be a complete Riemannian manifold equipped with a family of time dependent metrics $g(t)$ and potentials $\phi(t)$, $t\in [0, T]$. Suppose that 
$$\partial_t g=2h,$$
and for some function $\alpha_K: [0, T]\rightarrow \mathbb{R}$, it holds
\begin{eqnarray}
e^{-4Kt} (h+Ric_{m, n}(L)+Kg)-e^{-2Kt}h\geq  \alpha_K(t)g.
\end{eqnarray}
Suppose that
 $Ric_{m, n}(L)\geq -K_1$, $h\geq -K_2$, where $K_1, K_2$ are  two positive constants. Then
\begin{eqnarray*}
(L-\partial_t)\psi\geq -C_1K_2\psi^{1/2}-{C_1\over R}(m-1)\sqrt{K_1}\coth(\sqrt{K_1}\rho)-{C_2\over R^2}.
\end{eqnarray*}
\end{lemma}
{\it Proof}. By \cite{Li05}, under the condition $Ric_{m, n}(L)\geq -K_1$, on $(M, g(t), \phi(t))$,  it holds
$$
Ld(x_0, x, t)\leq (m-1)\sqrt{K_1}\rho\coth(\sqrt{K_1}\rho),
$$
and
\begin{eqnarray*}
L\psi&=&\eta'(d(x_0, x, t)/R){L d(x_0, x, t)\over R}+\eta''(d(x_0, x, t)/R){|\nabla d(x_0, x, t)|^2\over R^2}\\
&\geq &-{C_1\over R}(m-1)\sqrt{K_1}\coth(\sqrt{K_1}\rho)-{C_2\over R^2}.
\end{eqnarray*}
On the other hand, let $\gamma: [a, b]\rightarrow M$ be a fixed path such that $\gamma(a)=x$ and $\gamma(b)=y$. Let $S=\dot \gamma(s)$. Given a time $t_0\in [0, T]$, assuming that $\gamma$ is parameterize by the arc length with respect to metric $g(t_0)$ on $M$, then $|S|=1$ at time $t=t_0$. Moreover,  the evolution of the length of $\gamma$ with respect to $g(t)$ is given by
\begin{eqnarray*}
\left.{d\over dt}\right|_{t=t_0}L_{g(t)}(\gamma)&=&\int_a^b \left.{d\over dt}\right|_{t=t_0}\sqrt{g(t)(S, S)}ds\\
&=&{1\over 2}\int_a^b \left.{\partial_tg(t)(S, S)\over \sqrt{g(t)(S, S)}}\right|_{t=t_0}ds\\
&=&{1\over 2}\int_a^b \left.{\partial g(t)\over \partial t}(S, S)\right|_{t=t_0}ds.
\end{eqnarray*}
This yields, under the assumption $h\geq -K_2$, where $K_2\geq 0$,
\begin{eqnarray*}
\partial_t d(x, y, t)= \int_a^b h(S, S)ds\geq -K_2d(x, y, t).
\end{eqnarray*}
Since $-C_1\eta^{1/2}(r)\leq \eta'(r)\leq 0$, and $K_2\geq 0$, it holds
\begin{eqnarray*}
-\partial_t \psi&=&-{\eta'(\rho/R)\partial_t d(x_0, x, t)\over R}\\
&\geq&{\eta'(\rho/R) K_2 d(x_0, x, t)\over R}\\
&\geq &-{C_1K_2\over R}\psi^{1/2} d(x_0, x,t).
\end{eqnarray*}
Combining this with the lower bound of $L\psi$, we have
\begin{eqnarray*}
(L-\partial_t)\psi\geq -C_1K_2\psi^{1/2}-{C_1\over R}(m-1)\sqrt{K_1}\coth(\sqrt{K_1}\rho)-{C_2\over R^2}.
\end{eqnarray*}
The proof of Lemma \ref{comparison2} is completed. \hfill $\square$

\subsection{The Li-Yau Harnack inequality}

In this subsection we prove the Li-Yau Harnack inequality for
the positive solution to the heat equation  $\partial_t u=L_tu$ of
the time dependent Witten Laplacian on complete Riemannian manifolds
equipped with a backward $(\alpha, K, m)$-super Ricci flow.

Let $u$ be a positive solution to the heat equation $\partial_t
u=L_tu$. Let $f=\log u$. Then $(\partial_t - L)f = |\nabla f|^{2}$.
For any $\alpha>1$, let $F=t(|\nabla f|^2-\alpha f_t)$.

Since $\rho$ is Lipschitz on the complement of the cut locus of $o$,
$\psi$ is a Lipschitz function with support in $Q_{2R, T}$. As
explained in Li and Yau [43], an argument of Calabi \cite{Cal58}
allows us to apply the maximum principle to $\psi F$. Let $(x_0,
t_0)\in M\times [0, T]$ be a point where $\psi F$ achieves the
maximum. Then, at $(x_0, t_0)$,
\begin{eqnarray*}
\partial_t(\psi F)\geq 0, \ \Delta(\psi F)\leq 0, \ \nabla(\psi F)=0,
\end{eqnarray*}
which yields
\begin{eqnarray*}
(L-\partial_t)(\psi F)=\Delta (\psi F)-\nabla\phi\cdot\nabla(\psi
F)-\partial_t(\psi F)\leq 0.
\end{eqnarray*}Note that
\begin{eqnarray*}
(L-\partial_t)(\psi F)=\psi(L-\partial_t)F+(L-\partial_t)\psi
F+2\nabla\psi \cdot\nabla F.
\end{eqnarray*}
By Lemma \ref{comparison2}, we have
\begin{eqnarray*}
(L-\partial_t)\psi\geq -C_1K_2\psi^{1/2}-{C_1\over
R}(m-1)\sqrt{K_1}\coth(\sqrt{K_1}\rho)-{C_2\over R^2}.
\end{eqnarray*}
Therefore, at $(x_0, t_0)$, we have
\begin{eqnarray}
0\geq \psi(L-\partial_t)F+2\nabla\psi \cdot\nabla F-A(R,
T)F,\label{mmm}
\end{eqnarray}
where
$$A(R, T):=C_1K_2\psi^{1/2}+{C_1\over R}(m-1)\sqrt{K_1}\coth(\sqrt{K_1}\rho)+{C_2\over R^2}.$$
Denote $$C_3=C_3(m, K_1, K_2, R, T)=A(R,
T)+2|\nabla\psi|^2\psi^{-1}.$$ We have
\begin{eqnarray*}
C_3(m, K_1, K_2, R, T)\leq C_1K_2+{C\over R}+{C_2\over R^2}.
\end{eqnarray*}
Note that, at $(x_0, t_0)$, $\nabla\psi\cdot\nabla F=-\psi |\nabla
\psi|^2 F$. Substituting $(\ref{ooooo})$ into $(\ref{mmm})$, at
$(x_0, t_0)$, we have
\begin{eqnarray*}
0&\geq&\psi(L-\partial_t)F-A(R, T)F+2\nabla\psi\cdot\nabla F\\
&\geq &\psi(L-\partial_t)F -(A(R, T) + 2|\nabla\psi|^2\psi^{-1})F\\
&\geq& {2\psi F^2\over m\alpha^2 t}-\left({\psi \over
t}+C_3\right)F+{4(\alpha-1) \psi |\nabla f|^2 F\over
m\alpha^2}-{2C_2\over R}\psi^{1/2}|\nabla f|F\\
& &+\psi t\left[{2(\alpha-1)^2\over m\alpha^2}|\nabla f|^4-2K|\nabla
f|^2-\alpha B |\nabla f|-{\alpha^2 A^2\over 2}\right].
\end{eqnarray*}
By the inequality $ax^2-bx\geq {4b^2\over a}$ and $(\ref{0})$, and
multiplying the both sides by $\psi t_0$, we have
\begin{eqnarray*}
0&\geq& {2(\psi F)^2\over m\alpha^2}-\left(\psi +C_3t+{m\alpha^2 C_2^2 t\over 4(\alpha-1)R^2}\right)\psi F\\
& &-\psi^2 t^2\left({m\alpha^2 (2K+\gamma)^2\over
8(\alpha-1)^2}+{\alpha^2 B^2 \over 4\gamma}+{\alpha^2 A^2\over
2}\right).
\end{eqnarray*}
This yields that, for any $(x, t) \in Q_{R, T}$,
\begin{eqnarray*}
F(x,t) &\leq& (\psi F)(x_0, t_0)\\
&\leq&{m\alpha^2\over 2}\left[1+C_3t_0+{m\alpha^2 C_2^2 t_0\over
4(\alpha-1)R^2}+\sqrt{ {2\over m\alpha^2}\psi^2
t_0^2\left({m\alpha^2 (2K+\gamma)^2\over 8(\alpha-1)^2}+{\alpha^2
B^2 \over 4\gamma}+{\alpha^2 A^2\over 2}\right)}\right]\\
&\leq& {m\alpha^2\over 2}\left[1+\left(C_1K_2+{C\over R}+{C_2\over
R^2}+{m\alpha^2 C_2^2 \over
4(\alpha-1)R^2}\right)T+T\sqrt{{(2K+\gamma)^2\over
4(\alpha-1)^2}+{B^2 \over 2m\gamma}+{A^2\over m}} \right].
\end{eqnarray*}
Letting $R\rightarrow \infty$, we can derive, for all $\gamma>0$, we
have
\begin{eqnarray*}
F(x,t) \leq {m\alpha^2\over
2}\left[1+\left(C_1K_2+\sqrt{{(2K+\gamma)^2\over 4(\alpha-1)^2}+{B^2
\over 2m\gamma}+{A^2\over m}}\right)T \right].
\end{eqnarray*}
Therefore
\begin{eqnarray*}
{|\nabla u|^2\over u}-\alpha {\partial_t u\over u} \leq
{m\alpha^2\over 2t}\left[1+\left(C_1K_2+\sqrt{{(2K+\gamma)^2\over
4(\alpha-1)^2}+{B^2 \over 2m\gamma}+{A^2\over m}}\right)T \right].
\end{eqnarray*}

\subsection{The Li-Yau-Hamilton Harnack inequality}

In this subsection we prove the Li-Yau-Hamilton Harnack
inequality for the positive solution to the heat equation
$\partial_t u=L_tu$ of the time dependent Witten Laplacian on
complete Riemannian manifolds equipped with a variant of the $(K, m)$-super Ricci
flow.

Let $u$ be a positive solution to the heat equation $\partial_t
u=L_tu$. Let $f=\log u$. Then $(\partial_t - L)f = |\nabla f|^{2}$.
Let
\begin{eqnarray*}
F=te^{-4Kt}|\nabla f|^2-te^{-2Kt}f_t.
\end{eqnarray*}
Since $\rho$ is Lipschitz on the complement of the cut locus of $o$, $\psi$ is a Lipschitz function with support
in $Q_{2R, T}$. As explained in Li and Yau [43], an argument of Calabi \cite{Cal58} allows us
to apply the maximum principle to $\psi F$. Let $(x_0, t_0)\in M\times [0, T]$ be
a point where $\psi F$ achieves the maximum. Then, at $(x_0, t_0)$,
\begin{eqnarray*}
\partial_t(\psi F)\geq 0, \ \Delta(\psi F)\leq 0, \ \nabla(\psi F)=0,
\end{eqnarray*}
which yields
\begin{eqnarray*}
(L-\partial_t)(\psi F)=\Delta (\psi F)-\nabla\phi\cdot\nabla(\psi
F)-\partial_t(\psi F)\leq 0.
\end{eqnarray*}Note that
\begin{eqnarray*}
(L-\partial_t)(\psi F)=\psi(L-\partial_t)F+(L-\partial_t)\psi F+2\nabla\psi \cdot\nabla F.
\end{eqnarray*}
By Lemma \ref{comparison2}, we have
\begin{eqnarray*}
(L-\partial_t)\psi\geq -C_1K_2\psi^{1/2}-{C_1\over R}(m-1)\sqrt{K_1}\coth(\sqrt{K_1}\rho)-{C_2\over R^2}.
\end{eqnarray*}
Therefore, at $(x_0, t_0)$, we have
\begin{eqnarray*}
0\geq \psi(L-\partial_t)F+2\nabla\psi \cdot\nabla F-A(R, R)F,
\end{eqnarray*}
where
$$A(R, T):=C_1K_2\psi^{1/2}+{C_1\over R}(m-1)\sqrt{K_1}\coth(\sqrt{K_1}\rho)+{C_2\over R^2}.$$
Denote $$C(n, K, R, T)=A(R, T)+2|\nabla\psi|^2\psi^{-1}.$$
We have
\begin{eqnarray*}
C(m, K_1, K_2, R, T)\leq C_1K_2+{C\over R}+{C_2\over R^2}.
\end{eqnarray*}
Substituting $(\ref{GGGG1})$ into $(\ref{mmm})$, we have
\begin{eqnarray*}
0&\geq&\psi(L-\partial_t)F-A(R, T)F+2\nabla\psi\cdot\nabla F\\
&\geq &\psi(L-\partial_t)F -(A(R, T) + 2|\nabla\psi|^2\psi^{-1})F\\
&\geq &\psi\left[\frac{2[te^{-2Kt}(e^{-2Kt} -1)|\nabla f|^2 - F]^2}{mt}- 2\langle \nabla F, \nabla f \rangle +  \frac{(2Kt-1)}{t}F\right] -C(m, K_1, K_2, R, T)F\\
& &\ \ \ +\psi t\left[-{A^2\over 2}+2\alpha_K(t)|\nabla f|^2-e^{-2Kt}B|\nabla f|\right]\\
&\geq& \psi\left[{2F^2\over mt}+\frac{2Kt-1}{t}F - 2\langle \nabla F, \nabla f\rangle -{tA^2\over 2} +\frac{4e^{-2Kt}(1 - e^{-2Kt})|\nabla f|^2 }{m}F \right]-C(m, K_1, K_2, R, T)F\\
& &\ \ \ +\psi t \left[{2\over m}e^{-4Kt}(1-e^{-2Kt})^2|\nabla f|^4+2\alpha_K(t)|\nabla f|^2-e^{-2Kt}B|\nabla f|\right]\\
&\geq& \frac{2 \psi}{mt}F^2 + \frac{4e^{-2Kt}(1 - e^{-2Kt})|\nabla f|^2  \psi}{m}F + 2F\langle \nabla \psi, \nabla f \rangle +  \left[(2K -\frac{1}{t})\psi -C(m, K_1, k_2, R, T)\right]F\\
& &\ \ \ +\psi t\left[ {2\over m}e^{-4Kt}(1-e^{-2Kt})^2|\nabla f|^4+2\alpha_K(t)|\nabla f|^2-e^{-2Kt}B|\nabla f|-{A^2\over 2} \right]\\
&\geq& \frac{2 \psi}{mt}F^2 + \frac{4e^{-2Kt}(1 - e^{-2Kt})|\nabla f|^2 \psi}{m}F - 2F|\nabla \psi||\nabla f| +  \left[(2K -\frac{1}{t})\psi -C(m, K_1, K_2, R, T)\right]F\\
& &\ \ \ +\psi t\left[ {2\over m}e^{-4Kt}(1-e^{-2Kt})^2|\nabla f|^4+2\alpha_K(t)|\nabla f|^2-e^{-2Kt}B|\nabla f|-{A^2\over 2} \right]\\
&\geq& \frac{2 \psi}{mt}F^2 + \frac{4e^{-2Kt}(1 - e^{-2Kt})|\nabla f|^2 \psi}{m}F - 2{C_2 \over R}F\psi^{1/2}|\nabla f| +  \left[(2K -\frac{1}{t})\psi -C(m, K_1, K_2, R, T)\right]F\\
& &\ \ \ +\psi t\left[ {2\over m}e^{-4Kt}(1-e^{-2Kt})^2|\nabla f|^4+2\alpha_K(t)|\nabla f|^2-e^{-2Kt}B|\nabla f|-{A^2\over 2} \right].
\end{eqnarray*}
Multiplying by $t$ on both sides,
\begin{eqnarray*}
0 &\geq& \frac{2\psi}{m}F^2 + t\psi\frac{4e^{-2Kt}(1 - e^{-2Kt})|\nabla f|^2}{m}F - 2t{C_2 \over R}F\psi^{1/2}|\nabla f| +  [(2Kt - 1)\psi - C(m, K_1, K_2, R, T)t]F\\
& &\ \ \ +\psi t^2\left[ {2\over m}e^{-4Kt}(1-e^{-2Kt})^2|\nabla f|^4+2\alpha_K(t)|\nabla f|^2-e^{-2Kt}B|\nabla f|-{A^2\over 2} \right]\\
&=& \frac{2\psi}{m}F^2 + tF\left[\psi\frac{4e^{-2Kt}(1 - e^{-2Kt})|\nabla f|^2}{m} - 2{C_2 \over R}\psi^{1/2}|\nabla f|\right] +  [(2Kt - 1)\psi - C(m, K_1, K_2, R, T)t]F\\
& & \ \ \ +\psi t^2\left[ {2\over m}e^{-4Kt}(1-e^{-2Kt})^2|\nabla f|^4+2\alpha_K(t)|\nabla f|^2-e^{-2Kt}B|\nabla f|-{A^2\over 2} \right]\\
&\geq& \frac{2\psi}{m}F^2 +  [(2Kt - 1)\psi - C(m, K_1, K_2, R, T)t]F\\
& &+ tF\left[\psi\frac{4e^{-2Kt}(1 - e^{-2Kt})|\nabla f|^2}{m} - {C_2m \over 4e^{-2Kt}(1 - e^{-2Kt})R^2} - \frac{4e^{-2Kt}(1 - e^{-2Kt})}{m}\psi|\nabla f|^2\right] \\
 & & \ \ \ +\psi t^2\left[ {2\over m}e^{-4Kt}(1-e^{-2Kt})^2|\nabla f|^4+2\alpha_K(t)|\nabla f|^2-e^{-2Kt}B|\nabla f|-{A^2\over 2} \right]\\
&\geq& \frac{2\psi}{m}F^2 +  \left[(2Kt - 1)\psi - C(m, K_1, K_2, R, T)t - {C_2 m\over 4e^{-2Kt}(1 - e^{-2Kt})R^2}t\right]F\\
& &\ \ \  -\psi t^2\left[{m(2\alpha_K(t)-\gamma)^2\over 8e^{-4Kt}(1-e^{-2Kt})^2}+{e^{-4Kt}B^2\over 4\gamma}+{A^2\over 2} \right].
\end{eqnarray*}
Notice that the above calculation is done at the point $(x_0, t_0)$.
Since $\psi F$ reaches its maximum at this point, we can assume that
$\psi F(x_0, t_0) > 0$. Thus
\begin{eqnarray*}
0 &\geq& \frac{2}{m}(\psi F)^2  -\left[1+ C(n, K, R, T)t + {C_2 \over
4e^{-2Kt}(1 - e^{-2Kt})R^2}t\right](\psi F)\\
& &\ \ \  -\psi^2 t^2\left[{m(2\alpha_K(t)-\gamma)^2\over 8e^{-4Kt}(1-e^{-2Kt})^2}+{e^{-4Kt}B^2\over 4\gamma}+{A^2\over 2} \right].
\end{eqnarray*}
This yields that, for any $(x, t) \in Q_{R, T}$,
\begin{eqnarray*}
F(x,t) &\leq& (\psi F)(x_0, t_0)\\
& \leq& \frac{m}{2}\left[1 + C(m, K_1, K_2, R, T)t_0 + {C_2mt_0 \over 4e^{-2Kt_0}(1 - e^{-2Kt_0})R^2}\right]\\
& & +{\sqrt{m}\over 2} \sqrt{A^2t_0^2+{m(2\alpha_K(t_0)-\gamma)^2t_0^2\over 4e^{-4Kt_0}(1-e^{-2Kt_0})^2}+{t_0^2e^{-4Kt_0}B^2\over 2\gamma} }\\
&\leq& \frac{m}{2} \left[1 + C(m, K_1, K_2, R, T)T + {C_2mT \over 4e^{-2KT}(1 -e^{-2KT})R^2}\right]\\
& & +{\sqrt{m}\over 2} \sqrt{A^2T^2+\max\limits_{t\in [0, T]} \left( {m(2\alpha_K(t)-\gamma)^2t^2\over 4e^{-4Kt}(1-e^{-2Kt})^2}+{t^2e^{-4Kt}B^2\over 2\gamma}\right)   }.
\end{eqnarray*}
Let $R \rightarrow \infty$, we obtain
$$
F \leq \frac{m}{2}\left[1+C_1K_2T+\max\limits_{t\in [0, T]} {|2\alpha_K(t)-\gamma|t\over 2e^{-2Kt}(1-e^{-2Kt})}\right]+{\sqrt{m}\over 2}\left(A+{B\over \sqrt{2\gamma}}\right)T.
$$
In the case $\alpha_K(t)=0$, i.e., $e^{-4Kt} (h+Ric_{m, n}(L)+K)-e^{-2Kt}h\geq  0$, we have
$$
F \leq \frac{m}{2}\left[1+C_1K_2T+\max\limits_{t\in [0, T]} {\gamma t\over 2e^{-2Kt}(1-e^{-2Kt})}\right]+{\sqrt{m}\over 2}\left(A+{B\over \sqrt{2\gamma}}\right)T.
$$
In addition, when $B=0$, we can take $\gamma\rightarrow 0$, and we have
$$
F \leq \frac{m}{2}\left[1+C_1K_2T+{AT\over \sqrt{m}}\right].
$$

\medskip

We can also extend Hamilton's second order estimate  (i.e., Theorem \ref{HH3}, Theorem \ref{HHH} and Theorem \ref{HHH-H})  to positive solutions to the heat equation associated with time dependent Witten Laplacian on complete Riemannian manifolds with variant of the $(K, m)$-super Ricci flow. To save the length of the paper, we omit it here.

\medskip

\noindent{\bf Acknowledgement}.  The authors would like to
thank D. Bakry, J.-M. Bismut, D. Elworthy, M. Ledoux, N. Mok, K.-T. Sturm, A. Thalmaier and F.-Y. Wang for helpful discussions and warm encouragements during the past
years, and Dr. Y.-Z. Wang for careful preview
and useful comments on the earlier versions of this paper.
Part of this work has been done when the second author visited l'Institut des Hautes Etudes Scientifiques and l'Institut des Math\'ematiques de Toulouse de l'Universit\'e Paul Sabatier during November-December 2014, and l'Universit\'e Paris XIII
during January-February 2016. He would like to thank Prof. D. Bakry, Prof. J-M. Bismut and Prof. F. Nier for making this visit possible, and to thank IHES, UPS and Univ. Paris 13 for providing very nice environment  to finish this work.

\medskip

\begin{flushleft}
\medskip\noindent

Songzi Li, School of Mathematical Science, Fudan University, 220, Handan Road,
Shanghai, 200432, China, and Institut de Math\'ematiques, Universit\'e Paul Sabatier
118, route de Narbonne, 31062, Toulouse Cedex 9, France

\medskip

Xiang-Dong Li, {\sc  Academy of Mathematics and System Science, Chinese
Academy of Sciences, 55, Zhongguancun East Road, Beijing, 100190, P. R. China,}\\
E-mail: xdli@amt.ac.cn
\end{flushleft}

\end{document}